\documentclass[aop,MSNbibl,citesort,dvips]{arximspdf}
\usepackage{mathbh}

\doi{10.1214/09-AOP499}
\volume{38}
\issue{2}
\pubyear{2010}
\firstpage{927}
\lastpage{953}

\makeatletter

\newtheorem{theo}{Theorem}
\newtheorem{prop}[theo]{Proposition}

\newtheorem{lemma}[theo]{Lemma}

\newproclaim{Remark}{Remark}

\makeatother

\begin{document}
\begin{frontmatter}

\title{Large deviations for intersection local times~in~critical
dimension}
\runtitle{Intersection local times}

\begin{aug}
\author[A]{\fnms{Fabienne} \snm{Castell}\corref{}\ead[label=e1]{castell@cmi.univ-mrs.fr}}
\runauthor{F. Castell}
\affiliation{Universit\'{e} de Provence}
\address[A]{LATP\\
UMR CNRS 6632, C.M.I.\\
Universit\'{e} de Provence\\
39 Rue Joliot-Curie\\
F-13453 Marseille cedex 13\\
France\\
\printead{e1}} 
\end{aug}

\received{\smonth{1} \syear{2009}}
\revised{\smonth{8} \syear{2009}}

%
\begin{abstract}
Let $(X_t, t\geq0)$ be a continuous time simple random walk on
$\mathbb{Z}^d$ ($d \geq3$), and let $ l_T(x) $ be the time spent by
$(X_t, t\geq0)$ on the site $x$ up to time $T$. We prove a large
deviations principle for the $q$-fold self-intersection local time $
I_T = \sum_{x \in\mathbb{Z}^d} l_T(x)^q $ in the critical case
$q=\frac{d}{d-2}$. When $q$ is integer, we obtain similar results for
the intersection local times of $q$ independent simple random walks.
\end{abstract}

%

\end{frontmatter}

\section{Introduction}

\subsection*{Position of the problem}

Let $(X_t, t\geq0)$ be a continuous time simple random walk on
$\mathbb{Z}^d$, whose generator is denoted by $\triangle$ [where
$\triangle f (x) \stackrel{\triangle}{=}\sum_{y \sim x} (f(y) -
f(x))$]. Let
\[
l_T(x) = \int_0^T \delta_x(X_s) \,ds.
\]
The quantity of interest in this paper is the so called $q$-fold
self-intersection local time
\[
I_T = \sum_{x \in\mathbb{Z}^d} l_T(x)^q.
\]
When $q$ is integer, then
\[
I_T = q! \int_{0 \leq s_1 \leq\cdots\leq s_q \leq T}
\delta_{X_{s_1}=X_{s_2}=\cdots=X_{s_q}} \,ds_1 \cdots ds_q,
\]
which measures the amount of time the random walk spends on sites
visited at least $q$-times. Quantities measuring how much a random walk
does intersect itself, such as the range of the random walk, or the
self-intersection local time, appear in many models in physics. Far
from being exhaustive, we can cite the Polaron problem (see, for
instance, \cite{DVPol,mansmann}), models of polymers (see, for
instance, \cite{BRev,W1,W2,W3}), or models of diffusion in random environments
\cite{asselahcastell,AC,B,FC-FP,FC,KS,Kosh}. Partly motivated by the understanding
of these models, many studies have been devoted to such quantities for
more than twenty years. To describe the known results, we focus on
$I_T$ in the case $q=2$, where the literature is more complete, and we
refer the reader to the monograph \cite{chen-book} in preparation for a
very complete exposition of the subject, including results on the
range, or intersection local times of independent random walks.

Regarding the typical behavior of $I_T$ for large $T$, the results
depend of course on the dimension $d$, and of the transience/recurrence
of the random walk. They are summarized in Table \ref{table1}, where
$\gamma_1$ and $\underline{\gamma}_1$ are,\vspace*{1pt} respectively, the
intersection local time and renormalized intersection local time of the
Brownian motion up to time~1, and $\sigma(d)$ is a constant depending
on the dimension $d$:

%
\begin{table}
\caption{Typical behavior of $I_T$ for $q=2$}
\label{table1}
\begin{tabular*}{\tablewidth}{@{\extracolsep{\fill}}lccc@{}}
\hline
$\bolds d$ & \textbf{Order of} $\bolds{E(I_T)}$
& \textbf{Convergence in law} & \textbf{References}
\\
\hline
$d=1$ & $ T^{3/2}$ & $\frac{I_T}{T^{{3/2}}}
\stackrel{(d)}{\longrightarrow}
\gamma_1$ & \cite{borodin,brydges-slade,chen-li,perkins}
\\
$d=2$ & $T\log(T)$ & $\frac{I_T - E(I_T)}{T} \stackrel
{(d)}{\longrightarrow}
\underline{\gamma}_1$ & \cite{dynkin,legall1,legall2,rosen,stoll}
\\
$d \geq3$ & $T$ &
$\frac{I_T - E(I_T)}{\sqrt{\operatorname{var}(I_T)}}
\stackrel{(d)}{\longrightarrow} \mathcal{N}(0,1)$,
&\quad \cite{chen08,jain-orey,jain-pruitt} \\[5pt]
& & $\operatorname{var}(I_T) \sim\cases{
\sigma(3) T \log(T), & si $d=3$, \cr
\sigma(d) T, & si $d \geq4$,}
$ &\\
\hline
\end{tabular*}
\end{table}

Once we know the typical behavior, on can ask for untypical ones, that is,
for the large and moderate deviations for $I_T$. In many models, such
as the Polaron problem or polymers models, this is actually the
question of interest. The table below is an attempt to summarize the
results for $q=2$, achieved in recent years concerning this problem.

%
\begin{table}[b]
\caption{Large and moderate deviations results for $I_T$ for $q=2$}
\label{table}
\begin{tabular*}{\tablewidth}{@{\extracolsep{\fill}}lccc@{}}
\hline
$\bolds d$ & $\bolds{P[{I_T - E(I_T) \geq b^2_T}]}$
& \textbf{Value of} $\bolds{b_T}$ & \textbf{References}
\\
\hline
$d\leq2$
& $\exp({- \frac{2}{\kappa_c(2,d)^{8/d}} T^{({d-4})/{d}}
b_T^{4/d} })$
& $T^{2-{d}/{2}} \ll b^2_T \ll T^2$
& \cite{bass-chen,chen-li,mansmann,legall3,bass-chen-rosen}
\\[3pt]
$d=3$ & $\exp({- \frac{b_T^4}{2 \sigma(3) T \log(T)}})$
& $ \sqrt{T \log(T)} \ll b_T^2 \ll\sqrt{T \log(T)^{3/2}} $
& \cite{chen-book} \\
& $\exp({- \frac{2}{\kappa_c(2,d)^{8/d}} T^{({d-4})/{d}} b_T^{4/d}})$
& $\sqrt{T \log(T)^{3/2}} \ll b^2_T \ll T^2$
& \cite{asselah07,chen-book}
\\
[3pt]
$d=4$ & $ \exp({- \frac{b_T^4}{2 \sigma(4) T}}) $
& $\sqrt{T} \ll b^2_T \leq\sqrt{T \log\log T}$
& \cite{jain-pruitt2} \\
[3pt]
$d\geq5$ & $ \exp({- \frac{b_T^4}{2 \sigma(d) T}}) $
& $\sqrt{T} \ll b^2_T \leq\sqrt{T \log\log T}$
& \cite{jain-pruitt2} \\
& $\exp({- c(d) b_T})$
& $ T \leq b^2_T \ll T^2$
& \cite{asselah08,AC}
\\
\hline
\end{tabular*}
\end{table}

In Table \ref{table}, $\kappa_c(2,d)$ is the best constant $c$
in the Gagliardo--Nirenberg inequality:
\[
\forall d \leq3, \exists c \in\,]0,\infty[, \mbox{s.t. }
\forall f\dvtx \mathbb{R}^d \mapsto\mathbb{R},\qquad
\|{f}\|_4 \leq c \|{f}\|_2^{1-{d}/{4}} \|{\nabla f}\|_2^{
{d}/{4}},
\]
while $c(d)$ is an explicit constant related to discrete variational
inequalities.

So the picture is now almost complete, except for the dimensions $d\geq
4$. Note the coexistence of two different regimes in dimensions $d=3$
and $d\geq5$. The first one is an extension of the central limit
theorem describing the typical behavior, the second one corresponds to
the same pattern than in dimension $d \leq2$. To understand it, we give
some heuristics in the general case for $q$, where we want to control
$P[{I_T - E(I_T) \geq b_T^q}]$. For $I_T$ to be atypically high, one
possible strategy for the random walk is to remain during a time
$\tau\leq T$, in a box of size $R$. If $\tau\gg R^2$, this event has a
probability of order $\exp(-\tau/R^2)$. If $\tau\gg R^d$, one can
expect that on the box of size $R$, the local time $l_{\tau}(x)$ is now
of order $\tau /R^d$, so that $I_T$ has increased of an amount of order
$\tau^q/R^{d(q-1)}=b^q_T$. Hence, $\tau=b_T R^{d/q'}$ where $q'$ is the
conjugate exponent of $q$. Therefore, this strategy has a probability
of order $\exp(-b_T R^{d/q'-2})$. The best choice for $R$ is now the
choice that maximizes $\exp(-b_T R^{d/q'-2})$, under the constraint $ T
\geq\tau\gg R^{\max(2,d)}$.
\begin{itemize}
\item If $d < 2q'$ or equivalently $q < \frac{d}{(d-2)_+}$,
the bigger is $R$, the bigger is $\exp(-b_T\times\break R^{d/q'-2})$, so that
the best strategy for the random walk to make $I_T$ of order $b_T^q$,
is to remain
all the time $T$ in a ball of radius of order $(T/b_T)^{q'/d}$,
leading to the result of Table \ref{table} for
$d\leq2$ and the second regime in $d=3$.
\item If $d>2q'$, the smaller is $R$, the bigger is $\exp(-b_T
R^{d/q'-2})$, so that
the best strategy for the random walk to make $I_T$ of order $b_T^q$,
is now to remain
during a time $\tau$ of order $b_T$ in a ball of radius $R$ of order $1$,
leading to the second regime of Table \ref{table} in $d\geq5$.
\item The case $d=2q'$ is critical.
In that case $\exp(-b_T R^{d/q'-2})$ does not depend
on~$R$, so that whatever the order of $R$, $1 \leq R \ll\sqrt{T/b_T}$,
the strategy
consisting to remain a time $\tau= b_T R^2$ in a ball of size $R$
has a probability
of order $\exp(-b_T)$. The critical feature of $d=2q'$ is also reflected
in the fact that the Gagliardo--Nirenberg inequality appearing in
the results for $d < 2q'$, is now replaced by the Sobolev inequality.
For these reasons,
there is no result concerning the large and moderate deviations of
$I_T$ for $d=2q'$.
\end{itemize}

\subsection*{Main results}

This paper is a contribution to the large and very large deviations for
$I_T$ in the critical case $d=2q'$. By large deviations, we mean
deviations of the order of the mean $E(I_T)$, and by very large, we
mean deviations of order much larger than the order of the mean. When
$q$ is an integer (i.e., when $d=3$ and $q=3$, or when $d=4$ and
$q=2$), we obtain also similar results for the mutual intersection
$Q_T$ of $q$ independent random walks $(X^{(i)}_t; t \geq0, 1 \leq i
\leq q)$, defined by:
\[
Q_T = \sum_{x \in\mathbb{Z}^d} \prod_{i=1}^q l^{(i)}_T(x)
= \int_{0 \leq s_1, \ldots, s_q \leq T}
\delta_{X^{(1)}_{s_1}=X^{(2)}_{s_2}=\cdots=X^{(q)}_{s_q}} \,ds_1
\cdots ds_q,
\]
where $ l^{(i)}_T(x) = \int_0^T \delta_x(X^{(i)}_s) \,ds$. To state our
main results, we introduce some notation. For any function $f\dvtx
\mathbb{Z}^d \mapsto\mathbb{R}$, $\|{f}\|_p$ is the $l_p$ norm of $f$
[$\|{f}\|_p^p =\sum _{x \in\mathbb{Z}^d} |f|^p(x)$], and $\nabla f$ is
the discrete gradient of $f$ [for all $j \in\{{1, \ldots, d}\}$, for
all $x \in\mathbb{Z}^d$, $\nabla_j f(x)=f(x+e_j)-f(x)$].
\begin{prop}
\label{SILT4} For $d \geq3 $,
let $C_S(d) \in\,]0;+\infty[$
be the best constant in the discrete Sobolev's inequality
\[
\forall f \in l^{{2d}/({d-2})}(\mathbb{Z}^d)\qquad
\|{f}\|_{{2d}/({d-2})} \leq C_S(d) \|{\nabla f}\|_2.
\]

\begin{enumerate}
\item Exponential moments for $I_T$.\vspace*{1pt}

Let $d \geq3$, and let $q=\frac{d}{d-2}$.
%
\begin{eqnarray}
\label{ExpUBI}
\mbox{If } T^{1/q} &\ll& b_T,
\forall\theta\in\biggl[0; \frac{1}{C^2_{S}(d)}\biggr[\qquad
\limsup_{T \rightarrow\infty} \frac{1}{b_T} \log E [{\exp({\theta
I_T^{1/q}})}] = 0.
\\
\label{ExpLBI}
\mbox{If } b_T &\ll& T,
\forall\theta> \frac{1}{C^2_{S}(d)}\qquad
\liminf_{T \rightarrow\infty} \frac{1}{b_T} \log E [{\exp({\theta
I_T^{1/q}})}] = + \infty.
\end{eqnarray}
\item Exponential moments for $Q_T$.

Assume that $d=4$ and $q=2$, or $d=3$ and $q=3$.
%
\begin{eqnarray}
\label{ExpUBQ}
\mbox{If } T^{{1/q}} &\ll& b_T,
\forall\theta\in\biggl[0; \frac{q}{C^2_{S}(d)}\biggr[\qquad
\limsup_{T \rightarrow\infty} \frac{1}{b_T} \log E [{\exp({\theta
Q_T^{1/q}})}] = 0.
\\
\label{ExpLBQ}
\mbox{If } b_T &\ll& T,
\forall\theta> \frac{q}{C^2_{S}(d)}\qquad
\liminf_{T \rightarrow\infty} \frac{1}{b_T} \log E [{\exp({\theta
Q_T^{1/q}})}] = + \infty.
\end{eqnarray}
\end{enumerate}
\end{prop}

From Proposition \ref{SILT4}, it is straightforward to obtain very
large deviations upper bounds for $I_T$ and $Q_T$.
However, due to the degenerate form of the log-Laplace of $I_T^{1/q}$,
the corresponding lower bounds are not a direct consequence of
Proposition~\ref{SILT4}. These lower bounds are actually the main statement of the
following theorem.
\begin{theo}[(Very large deviations)]
\label{ILT4VLD}
\begin{enumerate}
\item Very large deviations for $I_T$.

Assume that $d \geq3$, $q=\frac{d}{d-2}$, and $T \gg b_T \gg T^{1/q}$.
%
\begin{equation}
\label{SILT4VLD.eq}
\lim_{T \rightarrow\infty} \frac{1}{b_T} \log P [{I_T \geq b_T^q}]
= - \frac{1}{C_S^2(d)}.
\end{equation}
\item Very large deviations for $Q_T$.

Assume that $d=4$ and $q=2$, or $d=3$ and $q=3$,
and that $ T \gg b_T \gg T^{1/q}$.
%
\begin{equation}
\label{MILT4VLD.eq}
\lim_{T \rightarrow\infty} \frac{1}{b_T} \log P [{Q_T \geq b_T^q}]
= - \frac{q}{C_S^2(d)}.
\end{equation}
\end{enumerate}
\end{theo}

Concerning the large deviations, our result is less precise since
the lower and upper bounds are different. To state it, we recall that
for $d \geq3$ and $q>1$,
$\lim_{T \rightarrow\infty} \frac{1}{T} E[{I_T}]$
exists in $\mathbb{R}^+$ [when $q$ is integer,
this limit is equal to $q! G_d(0)^{q-1}$,
where $G_d$ is the Green kernel of the
simple random walk on $\mathbb{Z}^d$].
\begin{theo}[(Large deviations for $I_T$)]
\label{SILT4LD}
Assume that $d \geq3$, $q=\frac{d}{d-2}$. There exists a constant
$c(d) > 0$ such that $\forall y > c(d)$
\begin{eqnarray}
\label{SILT4LD.eq}
- \frac{y^{1/q}}{C_S^2(d)}
& \leq&
\liminf_{T \rightarrow\infty} \frac{1}{T^{1/q}} \log P[{I_T \geq T y}]
\nonumber\\[-8pt]\\[-8pt]
& \leq&
\limsup_{T \rightarrow\infty} \frac{1}{T^{1/q}} \log P[{I_T \geq T y}]
= - \frac{1}{c(d)} y^{1/q}.\nonumber
\end{eqnarray}
\end{theo}
\begin{Remark}\label{Remark1}
Unfortunately, our proof does not allow to obtain the result for all $y
> \lim_{T \rightarrow\infty} \frac{E(I_T)}{T}$.
\end{Remark}
\begin{Remark}\label{Remark2}
As in Theorem \ref{ILT4VLD}, we could obtain similar results for $Q_T$.
However, such a result would not correspond to a large deviations
result for $Q_T$, since $E(Q_T)$ is of order $\log(T)$ for $d \geq3$
and $q= d/(d-2)$. Concerning $Q_T$, we should also mention that papers
\cite{marcus-rosen-97} and \cite{rosen-97} give moderate deviations
estimates $P[{Q_T - E(Q_T) \geq\log(T) b_T}]$ for scales $b_T$ up to
$\log\log\log(T)$.
\end{Remark}

\subsection*{Sketch of the proof}

The proof of the lower bounds is easy and
relies heavily on the large deviations results for $\frac{l_T}{T}$ proved
by Donsker and Varadhan. Namely,
let $\mathcal{F}=\{{\mu\dvtx\mathbb{Z}^d \mapsto\mathbb{R}^+ ; \sum
_{x \in\mathbb{Z}^d} \mu(x)=1}\}$.
$\mathcal{F}$ is endowed with the weak topology of probability
measures. By
the results of Donsker and Varadhan \cite{DV76}, $l_T/T$ satisfy
a restricted large deviations principle in $\mathcal{F}$ (by
``restricted,'' it
is meant that the large deviations upper bound is only true for
compact sets), with rate function $\mathcal{I}(\mu)= \|{\nabla\sqrt
{\mu}}\|_2^2$.
Now, for any $M$ satisfying $Mb_T \leq T$,
${\frac{I_T}{b_T^q}} \geq{\frac{I_{Mb_T}}{b_T^q}}=M^q
\|{{\frac{l_{Mb_T}}{Mb_T}}}\|_q^q$.
Moreover, the function $\mu\in\mathcal{F}\mapsto\|{\mu}\|_q
=\sup\{{\sum_x \mu(x) f(x) ; f \mbox{ compactly supported, } \|
{f}\|_{q'}=1}\}$
is lower semicontinuous in weak topology. The large deviations lower bound
for $\frac{l_{Mb_T}}{Mb_T}$ [with the change of variable $\mu
(x)=g^2(x)$], leads
therefore to
\begin{eqnarray}
\label{sketchLB}
&&\liminf_{T \rightarrow\infty} \frac{1}{b_T} \log P [{I_T > b_T^q}]
\nonumber\\[-8pt]\\[-8pt]
&&\qquad\geq- M \inf\biggl\{{\|{\nabla g}\|_2^2 ; g \mbox{ such that } \|{g}\|_2
=1 \mbox{ and } \|{g}\|_{2q}^{2} > \frac{1}{M}}\biggr\}\nonumber
\end{eqnarray}
for all $M < \liminf\frac{T}{b_T}$. For $b_T \ll T$, all the values of
$M$ are
allowed, and taking the supremum in $M$ in (\ref{sketchLB})
leads to the lower bound in (\ref{SILT4VLD.eq}). Actually, this argument
remains valid for any scale $b_T$ such that $1 \ll b_T \ll T$ (see Proposition
\ref{SILT4LB}).

For the very
large deviations upper bound for $I_T$, the results of Donsker and
Varadhan are not sufficient, since on one hand, the large deviations upper
bound for $l_T/T$ is only true for compact sets of $\mathcal{F}$, and
on the
other hand,
the function $\mu\in\mathcal{F}\mapsto\|{\mu}\|_q$ is not continuous.
We present now the
main ingredients of the proof of the upper bound (\ref{ExpUBI}).
First of all, it is easy to see that $I_T \leq I_T(R)$,
the intersection local time
of the random walk folded on the torus of radius $R$. Now,
the main tool in the proof is
the mysterious Dynkin isomorphism theorem, according to which
the law of the local times
of a symmetric recurrent Markov process stopped at an independent
exponential time,
is related to the law of the square of a Gaussian
process whose covariance function is the Green kernel of the stopped
Markov process.
This allows us
to control the exponential moments of $I_T^{1/q}$, with
the exponential moments of $N_T(R)=\frac{1}{2}
({\sum_{x \in\mathbb{T}_R} Z_x^{2q}})^{1/q}=\frac{1}{2}
\|{Z}\|_{2q,R}^2$
where:
\begin{itemize}[-]
\item[-] $\mathbb{T}_R$ is the torus of radius $R$;
\item[-] $(Z_x, x \in\mathbb{T}_R)$ is a centered Gaussian process
whose covariance function is given by $G_{R,\lambda}(x,y)$,
the Green kernel of the simple
random walk on $\mathbb{T}_R$,
stopped at an independent exponential time
with parameter $\lambda\sim b_T/T$, (Lemmas
\ref{redtorus.lem}, \ref{dynkin.lem} and~\ref{ubl4norm.lem});
\item[-] $\|\cdot\|_{2q,R}$ denotes the norm in $l^{2q}(\mathbb{T}_R)$.
\end{itemize}

We can now rely on concentration inequalities for norms
of Gaussian processes.
Let $M_{R,T}$ denote the median of $\|{Z}\|_{2q,R}$. For small $\alpha$,
\[
\exp\biggl[{\frac{\theta}{2} \|{Z}\|_{2q,R}^2}\biggr]
\leq\exp\biggl[{\frac{\theta(1+\alpha)}{2} (\|{Z}\|_{2q,R}-M_{R,T})^2}\biggr]
\exp\biggl[{\frac{\theta(1+\alpha)}{2 \alpha} M_{R,T}^2}\biggr].
\]
By concentration inequalities, the tail behavior of $\|{Z}\|
_{2q,R}-M_{R,T}$ is
that of a centered Gaussian variable with variance
\[
\rho= \sup\{{\langle{f,G_{R,\lambda}f}\rangle ; \|{f}\|
_{(2q)',R} =1}\}
.
\]
Therefore, for $\theta< \frac{1}{(1+\alpha)\rho}$,
\[
\exp\biggl[{\frac{\theta(1+\alpha)}{2} (\|{Z}\|_{2q,R}-M{R,T})^2}\biggr]
\leq\frac{1}{\sqrt{1-\theta(1+\alpha) \rho}}.
\]
Besides,\vspace*{1pt} one can prove that $M_{R,T}$ is of order $R^{{d}/({2q})}$ as
soon as
$\lambda R^d \gg1$, and that $\rho\sim\frac{1}{C_S^2(d)}$ if
$\lambda R^2 \gg1$.
We therefore\vspace*{-2pt} obtain the result in (\ref{ExpUBI}), if $R$ is chosen so that
$b_T \gg R^{d/q}$ and $\lambda R^2 \sim\frac{b_T}{T} R^2 \gg1$. The
best choice
for $R$ is now to take $R^{d/q} = T/R^2$, i.e., $R=T^{1/d}$ since
$q=\frac{d}{d-2}$,
leading to $b_T \gg T^{1/q}$.

\subsection*{An open question}
The large, very large and moderate deviations for $I_T$ and $Q_T$
in the subcritical case (i.e., $d \leq2$, or $d=3$ and $q< \frac{d}{d-2}$)
are linked to Gagliardo--Nirenberg inequality
in a continuous setting (i.e., for functions $f$ from $\mathbb{R}^d$
to $\mathbb{R}$),
while the same problem in supercritical case $
d \geq3$ and $q > \frac{d}{d-2}$,
is linked to functional inequality in a discrete setting. One can therefore
think that in the critical case $q=\frac{d}{d-2}$, the moderate deviations
of $I_T - E[{I_T}]$ are at least up to some scale, related to
the Sobolev inequality in a
continuous setting. However, since the best constants in the discrete
and continuous Sobolev inequality are the same, this would not change
the statement. Therefore, we do believe that in the critical case
$d=2q'$, there are only two regimes of deviations from the mean:
\[
P[{I_T - E(I_T) \geq b_T^q}] \asymp\cases{
\exp\biggl({-\dfrac{b_T^{2q}}{2\sigma(d)T}}\biggr), &\quad for
$\sqrt{T} \ll b_T^q \ll T^{{q}/({2q-1})}$, \vspace*{2pt}\cr
\exp\biggl({-\dfrac{1}{C_S^2(d)} b_T}\biggr), &\quad for
$T^{{q}/({2q-1})} \ll b_T^q \ll T^q$.}
\]
We do not know how to prove this result. Actually, the same question is
also open in the supercritical case (with $\frac{1}{C_S^2(d)}$ replaced
by the constant $c(d)$ given in \cite{asselah08}).

The paper is organized as follows. Section \ref{LB} is devoted
to the proof of exponential moments lower bounds (\ref{ExpLBI}) and
(\ref{ExpLBQ}).
In Section \ref{UB}, we prove the exponential moments upper bounds
(\ref{ExpUBI}) and
(\ref{ExpUBQ}). In Section \ref{LBGD}, we give the proof of the large
and very
large deviations lower bounds. With Proposition \ref{SILT4}, this ends the
proof of Theorem \ref{ILT4VLD}. Finally, Section \ref{UBGD} is devoted
to the proof of the upper bound in (\ref{SILT4LD.eq}), which ends the proof
of Theorem \ref{SILT4LD}.


\section{Exponential moments lower bound}
\label{LB}

This section is devoted to the proof of the lower bounds (\ref{ExpLBI})
and (\ref{ExpLBQ})
in Proposition
\ref{SILT4}.

\subsection*{Lower bound for $I_T$}
Fix $M>0$. Since $b_T \ll T$, for $T$ sufficiently large [$T \geq
T_0(M)$]
$Mb_T \leq T$, and $I_T \geq I_{Mb_T}$.
For any $f$ such that $\|{f}\|_{q'}=1$,
%
\begin{equation}
E [{\exp(\theta I_T^{1/q})}]
\geq E [{\exp(\theta I_{Mb_T}^{1/q})}]
\geq E \biggl[{\exp\biggl({\theta\sum_x f(x)
l_{Mb_T}(x)}\biggr)}\biggr].
\end{equation}

It is a standard result that the occupation measure of $X$ satisfies a weak
large deviations principle in $\mathcal{F}$, in $\tau$-topology
(i.e., the
topology defined
by duality with bounded measurable functions), with rate function
$\mathcal{J}(\mu) = \|{\nabla\sqrt{\mu}}\|^2$ (see, for instance, Theorem
5.3.10, page 210 in \cite{DS}). Since $f$ is bounded by 1 as soon as
$\|{f}\|_{q'}=1$, the function $\mu\in\mathcal{F}\mapsto\sum_{x
\in\mathbb{Z}^d}
f(x) \mu(x)$ is continuous in
$\tau$-topology and the large deviations lower bound for\break $\frac{1}{Mb_T}
\int_0^{Mb_T} \delta_{X_s} \,ds$ (written with the change of variable
$g=\sqrt{\mu}$)
yields: $\forall\theta\geq0$,
$\forall M >0$,
$\forall f \in l_{q'}(\mathbb{Z}^d)$ such that $\|{f}\|_{q'}=1$,
%
\begin{equation}\qquad
\label{LBQ.eq}
\liminf_{T\rightarrow\infty} \frac{1}{b_T}
\log E [{\exp(\theta I_T^{1/q})}]
\geq M \sup_{g, \|{g}\|_2=1} \biggl\{{\theta\sum_x f(x) g^2(x) - \|{\nabla
g}\|_2^2}\biggr\}.
\end{equation}
Assume now that $\theta> \frac{1}{C_S^2(d)}
= \inf\frac{\|{\nabla f}\|_2^2}{\|{f}\|_{2q}^2}$ for $q=\frac{d}{d-2}$.
Since the infimum can be
reduced to the infimum over compactly supported functions $f$, we can
find $g_0$ with compact support in $\mathbb{Z}^d$, such that
$\theta> \frac{\|{\nabla g_0}\|_2^2}{\|{g_0}\|_{2q}^2}$. Dividing
$g_0$ by its $l_2$-norm if necessary, we can moreover assume that
$\|{g_0}\|_2=1$. We now take $f= \frac{g_0^{2(q-1)}}{\|{g_0}\|^{2(q-1)}_{2q}}$
(note that $\|{f}\|_{q'} =1$), $g=g_0$ in
(\ref{LBQ.eq}). $\forall M >0$,
\begin{eqnarray*}
\liminf_{T\rightarrow\infty} \frac{1}{b_T}
\log E [{\exp(\theta I_T^{1/q})}]
& \geq& M \biggl({\theta\sum_x f(x) g^2(x) - \|{\nabla g}\|_2^2}\biggr)
\\
& = & M \biggl({\theta\frac{\sum_x g_0^{2q}(x)}{\|{g_0}\|^{2(q-1)}_{2q}} -
\|{\nabla g_0}\|_2^2 }\biggr)
\\
& = & M ({\theta\|{g_0}\|_{2q}^2 - \|{\nabla g_0}\|_2^2 })
.
\end{eqnarray*}
But $\theta\|{g_0}\|_{2q}^2 - \|{\nabla g_0}\|_2^2 > 0$, so that
(\ref{ExpLBI}) is proved by sending $M$ to infinity.

\subsection*{Lower bound for $Q_T$}
Fix $M>0$. Since $b_T \ll T$, for $T$ sufficiently large [$T \geq
T_0(M)$]
$Mb_T \leq T$, and $Q_T \geq Q_{Mb_T}$.
$\forall\theta\geq0$, and $\forall m \in\mathbb{N}$,
\begin{eqnarray*}
E [{\exp(\theta Q_T^{1/q})}]
& \geq& E [{\exp(\theta Q_{Mb_T}^{1/q})}]
\\
& \geq& \frac{\theta^{qm}}{(qm)!} E [{Q_{Mb_T}^m}]
\\
& = & \frac{\theta^{qm}}{(qm)!} \sum_{x_1, \ldots, x_m}
E\Biggl[{\prod_{j=1}^q \prod_{i=1}^m l^{(j)}_{Mb_T}(x_i) }\Biggr]
\\
& = & \frac{\theta^{qm}}{(qm)!} \sum_{x_1, \ldots, x_m}
E\Biggl[{\prod_{i=1}^m l_{Mb_T}(x_i)}\Biggr]^q
\\
& \geq& \frac{\theta^{qm}}{(qm)!} \Biggl[{\sum_{x_1, \ldots, x_m} f(x_1)
\cdots f(x_m) E\Biggl[{\prod_{i=1}^m l_{Mb_T}(x_i)}\Biggr]}\Biggr]^q
\end{eqnarray*}
for any $f \in l_{q'}(\mathbb{Z}^d)$, such that $\|{f}\|_{q'}=1$.
Therefore, $\forall\theta\geq0$, and $\forall m \in\mathbb{N}$,
%
\begin{equation}
E [{\exp(\theta Q_T^{1/q})}]^{1/q}
\geq\frac{\theta^m}{((qm)!)^{1/q}} E\biggl[{\biggl({\sum_x f(x)
l_{Mb_T}(x)}\biggr)^m}\biggr].
\end{equation}
It follows from Stirling's formula that there exists $C>0$ such that
$\forall m \in\mathbb{N}$, $\frac{1}{((qm)!)^{1/q}} \geq C \frac
{1}{q^m m!}$.
Hence, $\forall\theta\geq0$, and $\forall m \in\mathbb{N}$,
%
\begin{equation}
E [{\exp(\theta Q_T^{1/q})}]^{1/q}
\geq C \frac{1}{m!} E\biggl[{\biggl({\frac{\theta}{q} \int_0^{Mb_T} f(X_s)
\,ds}\biggr)^m}\biggr].
\end{equation}
Summing over $m$, we have thus proved that for $T \geq T_0(M)$,
$\forall\theta\geq0$, $\forall f \in l_{q'}(\mathbb{Z}^d)$ such that
$\|{f}\|_{q'}=1$,
\[
E [{\exp(\theta Q_T^{1/q})}]^{1/q}
\geq C E \biggl[{\exp\biggl({\frac{\theta}{q} \int_0^{Mb_T} f(X_s) \,ds}\biggr)}\biggr].
\]
At this point, the proof is the same as the proof of the lower bound
for $I_T$.


\section{Exponential moments upper bounds}
\label{UB}
In this section, we obtain an upper bound for the exponential moments
of $I_T^{1/q}$ and $Q_T^{1/q}$.

\textit{Step} 1.
Comparison with the SILT of the random walk on the torus,
stopped at an exponential
time.
\begin{lemma}
\label{redtorus.lem}
Let $\alpha> 0$, and let $\tau$ be an exponential random variable
with parameter $\lambda=
\alpha\frac{b_T}{T}$, independent of the random walk $(X_s, s \geq0)$.
Let $R \in\mathbb{N}^*$, and
let us denote by $X^{(R)}_s = X_s \operatorname{mod}(R)$ the simple random walk on
$\mathbb{T}_R$,
the $d$-di\-men\-sional discrete torus of radius $R$. Finally,
let $l^{(R)}_{\tau}(x)=\int_0^{\tau}
\delta_x (X_s^{(R)}) \,ds$, and $I_{R,\tau}= \sum_{x \in\mathbb{T}_R}
(l^{(R)}_{\tau}(x))^q$.
Then, $\forall\theta> 0$, $\forall\alpha> 0$, $\forall R > 0$,
$\forall T >0$,
%
\begin{equation}
\label{redtorus.eq}
E [{\exp({\theta I_T^{1/q}})}]
\leq e^{\alpha b_T} E [{\exp({\theta I_{R,\tau}^{1/q}})}].
\end{equation}
\end{lemma}
\begin{pf}
\begin{eqnarray*}
I_T & = & \sum_{x \in\mathbb{Z}^d} l_T^q(x)
= \sum_{x \in\mathbb{T}_R} \sum_{k \in\mathbb{Z}^d} l^q_T(x+kR)
\\
& \leq& \sum_{x \in\mathbb{T}_R} \biggl({\sum_{k \in\mathbb{Z}^d} l_T(x+kR)}\biggr)^q
= \sum_{x \in\mathbb{T}_R} l^q_{R,T}(x) = I_{R,T}.
\end{eqnarray*}
Therefore,
\begin{eqnarray*}
E[{\exp({\theta I_T^{1/q}})}] \exp({-\alpha b_T})
& \leq&
E[{\exp({\theta I_{R,T}^{1/q}})}] P[{\tau\geq T}]
\\
& \leq&
E[{\exp({\theta I_{R,T}^{1/q}}) \mathbh{1}_{\tau\geq T}}]
\\
& \leq&
E[{\exp({\theta I_{R,\tau}^{1/q}})}],
\end{eqnarray*}
where the first inequality comes from the choice of
$\lambda=
\alpha\frac{b_T}{T}$, and the second one from independence
of $\tau$ and $X$.
\end{pf}

\textit{Step} 2. \textit{The Eisenbaum isomorphism theorem}.
There are various versions of isomorphism theorems
in the spirit of the Dynkin isomorphism theorem. We use here the following
version due to Eisenbaum \cite{eisenbaum} (see also Corollary 8.1.2,
page~364
in \cite{marcus-rosen}).
\begin{theo}[(Eisenbaum)]
\label{dynkin.lem}
Let $\alpha$ and $\tau$ be as in Lemma \ref{redtorus.lem}.
Let us define for all $x, y \in\mathbb{T}_R$,
$G_{R,\lambda}(x,y)=E_x [{\int_0^\tau\delta_y(X^{(R)}_s) \,ds}]$.
Let $(Z_x, x \in\mathbb{T}_R)$ be a centered Gaussian process
with covariance matrix $G_{R,\lambda}$,
independent of $\tau$ and of the random walk $(X_s, s \geq0)$.
For $s\neq0$, consider the process
$S_x := l^{(R)}_{\tau}(x) + \frac{1}{2} (Z_x+s)^2$. Then, for all measurable
and bounded function $F \dvtx \mathbb{R}^{\mathbb{T}_R} \mapsto\mathbb{R}$,
%
\begin{equation}
\label{dynkin.eq}
E\bigl[{F\bigl((S_x; x\in\mathbb{T}_R)\bigr)}\bigr]
= E\biggl[{F\biggl({\biggl(\frac{1}{2}(Z_x +s)^2;x \in\mathbb{T}_R\biggr)}\biggr) \biggl({1 + \frac
{Z_0}{s}}\biggr)}\biggr].
\end{equation}
\end{theo}

\textit{Step} 3.
Comparison between exponential moments of $I_T$ and
exponential moments for $\sum_x Z_x^{2q}$.

Theorem \ref{dynkin.lem} allows one to control exponential moments
of $I_{R,\tau}^{1/q}$ by exponential moments of
$(\sum_{x \in\mathbb{T}_R} Z_x^{2q})^{1/q}$.
\begin{lemma}
\label{ubl4norm.lem}
For any $\alpha>0$ and $R>0$, let $\tau$ and $(Z_x, x\in\mathbb
{T}_R)$ be defined
as in Lemma \ref{dynkin.lem}.
$\forall\alpha> 0$, $\forall\theta>0$, $\forall\gamma> \theta$,
$\forall\varepsilon\in\,]0; \min(1, \sqrt{\frac{\gamma}{\theta} -1})[$,
$\forall R > 0$, $\forall T >0$, there exists a constant $C(\varepsilon)
\in\,]0;\infty[$
depending only on $\varepsilon$, such that
%
\begin{eqnarray}
\label{ubl4norm.eq}
&&E [{\exp({\theta I_{R,\tau}^{1/q}})}]
\nonumber\\
&&\qquad\leq
1 + C(\varepsilon) {\frac{\theta}{\gamma- \theta(1+\varepsilon)^2}}
\biggl({1+{\frac{\sqrt{T } R^{d/2q}}{\sqrt{\alpha} b_T}}}\biggr)
\\
&&\hspace*{48.1pt}{}\times
{\frac{E[{\exp({{{\gamma}/{2}} \|{Z}\|_{2q,R}^2})}]^{
{1}/({1+\varepsilon})}}{P[{\|{Z}\|_{2q,R} \geq2 \sqrt{2b_T\varepsilon}}]}}
\exp({\gamma\varepsilon^2 b_T}),\nonumber
\end{eqnarray}
where $\|\cdot\|_{p,R}$ is the
$l_p$ norm of functions on $\mathbb{T}_R$.
\end{lemma}
\begin{pf}
By independence of $(Z_x, x \in\mathbb{T}_R)$ and $(X_s, s \geq0)$,
$\forall s\neq0$, $\forall y >0$, $\forall\varepsilon> 0$,
%
\begin{eqnarray}
\label{UB_I.eq}\quad
&&P\biggl[{\sum_{x \in\mathbb{T}_R} {\frac{(Z_x+s)^{2q}}{2^q}} \geq b_T^q
\varepsilon^q}\biggr]
P[{I_{R,\tau} \geq b_T^q y^q}]
\nonumber\\
&&\qquad= P\biggl[{\sum_{x \in\mathbb{T}_R} {\frac{(Z_x+s)^{2q}}{2^q}} \geq b_T^q
\varepsilon^q ; \sum_{x \in\mathbb{T}_R} \bigl(l^{(R)}_{\tau
}(x)\bigr)^q \geq b_T^q y^q}\biggr]
\nonumber\\[-8pt]\\[-8pt]
&&\qquad\leq P\biggl[{\sum_{x \in\mathbb{T}_R} S_x^q \geq b_T^q (y^q+\varepsilon^q)}\biggr]
\nonumber\\
&&\qquad= E\biggl[{\biggl(1+ {\frac{Z_0}{s}}\biggr) \mathbh{1}_{ \sum_{x \in\mathbb{T}_R}
{(Z_x+s)^{2q}}/{2^q} \geq b_T^q (y^q+\varepsilon^q)} }\biggr]\qquad
\mbox{by Theorem \ref{dynkin.lem}.}\nonumber
\end{eqnarray}
Hence, using Markov inequality,
%
\begin{eqnarray}
\label{ubexp.eq}
E[{\exp({\theta I_{R,\tau}^{1/q}})}]
&=& 1 + \int_0^{\infty} \theta b_T e^{\theta b_T y} P[{I_{R,\tau}
\geq b_T^q y^q}] \,dy
\nonumber\\
&\leq&
1 + {\frac{E [{({1+ {{Z_0}/{s}}}) \exp({{{\gamma}/{2}} \|
{Z +s\mathbh{1}}\|_{2q,R}^2})}]}{P[{\|{Z+s\mathbh{1}}\|_{2q,R} \geq
\sqrt{2 b_T \varepsilon}}]}}
\\
&&\hspace*{15.7pt}{}\times
\int_0^{\infty} \theta b_T e^{\theta b_T y} e^{-b_T \gamma
(y^q+\varepsilon
^q)^{1/q}}\,
dy.\nonumber
\end{eqnarray}

Now, $\forall\varepsilon> 0$, $\forall\theta>0$, $\forall\gamma>
\theta$,
$\forall T > 0$,
%
\begin{equation}
\label{ubintegrale.eq}\quad
\int_0^{\infty} \theta b_T
e^{\theta b_T y} e^{-b_T \gamma(y^q+\varepsilon^q)^{1/q}} dy
\leq\int_0^{\infty} \theta b_T e^{\theta b_T y} e^{-b_T \gamma y} dy
= {\frac{\theta}{\gamma- \theta}}.
\end{equation}
Regarding the denominator in (\ref{ubexp.eq}),
%
\begin{eqnarray}
P\bigl[{\|{Z +s \mathbh{1}}\|_{2q,R} \geq\sqrt{2b_T \varepsilon}}\bigr]
& \geq& P \bigl[{\|{Z}\|_{2q,R} \geq\sqrt{2b_T \varepsilon} + \|{s \mathbh
{1}}\|_{2q,R}}\bigr]
\\
\label{lbden.eq}
& = & P \bigl[{\|{Z}\|_{2q,R} \geq\sqrt{2b_T \varepsilon} + |s| R^{d/2q}}\bigr].
\end{eqnarray}

On the other hand, $\forall\varepsilon> 0$,
\[
\|{Z +s\mathbh{1}}\|_{2q,R}^2 \leq(\|{Z}\|_{2q,R} + \|{s\mathbh{1}}\|
_{2q,R})^2
\leq\|{Z}\|_{2q,R}^2 (1+\varepsilon) + \biggl(1 + {\frac{1}{\varepsilon}}\biggr)
\|{s\mathbh{1}}\|_{2q,R}^2,
\]
so that
%
\begin{eqnarray}
\label{ubnum.eq}
&& E \biggl[{\biggl({1+ {\frac{Z_0}{s}}}\biggr) \exp\biggl({{\frac{\gamma}{2}} \|{Z +s \mathbh
{1}}\|_{2q,R}^2}\biggr)}\biggr]
\nonumber\\
&&\qquad\leq E \biggl[{\biggl({1+ {\frac{Z_0}{s}}}\biggr) \exp\biggl({{\frac{\gamma}{2}}
(1+\varepsilon) \|{Z}\|_{2q,R}^2}\biggr)}\biggr] \exp\biggl({{\frac{\gamma}{2}} {\frac
{1+\varepsilon}{\varepsilon}} s^2 R^{d/q}}\biggr)
\nonumber\\[-8pt]\\[-8pt]
&&\qquad\leq
E \biggl[{\biggl|{1+ {\frac{Z_0}{s}}}\biggr|^{({1+\varepsilon})/{\varepsilon}}}\biggr]
^{{\varepsilon}/({1+\varepsilon})}
E\biggl[{\exp\biggl({{\frac{\gamma}{2}} (1+\varepsilon)^2 \|{Z}\|
_{2q,R}^2}\biggr)}\biggr]^{{1}/({1+\varepsilon})}\nonumber\\
&&\qquad\quad{}\times\exp\biggl({{\frac{\gamma}{2}} {\frac{1+\varepsilon}{\varepsilon}} s^2
R^{d/q}}\biggr),\nonumber
\end{eqnarray}
$Z_0$ being a centered Gaussian variable with variance $G_{R,\lambda}(0,0)
\leq E(\tau) = 1/\lambda$, for all $\varepsilon> 0$, there exists a constant
$C(\varepsilon)$ depending only on $\varepsilon$ such that
%
\begin{equation}
\label{ubdens.eq}
E \biggl[{\biggl|{1+ {\frac{Z_0}{s}}}\biggr|^{({1+\varepsilon})/{\varepsilon}}}\biggr]^{
{\varepsilon}/({1+\varepsilon})}
\leq C(\varepsilon) \Biggl({1+ \sqrt{{\frac{T}{\alpha b_T}}} {\frac{1}{s}}}\Biggr).
\end{equation}
Putting (\ref{ubexp.eq}), (\ref{ubintegrale.eq}), (\ref{lbden.eq}),
(\ref{ubnum.eq}) and (\ref{ubdens.eq}) together, we have thus proved that
$\forall\theta> 0$, $\forall\gamma> \theta$, $\forall\varepsilon>0$,
$\forall R > 0$, $\forall T > 0$, $\forall s \neq0$,
%
\begin{eqnarray}\quad
\label{ubexp2.eq}
&& E[{\exp(\theta I_{R,\tau}^{1/q})}]
\nonumber\\
&&\qquad \leq1 + C(\varepsilon) {\frac{\theta}{\gamma- \theta}}
\Biggl({1+ \sqrt{{\frac{T}{\alpha b_T}}}
{\frac{1}{s}}}\Biggr)\nonumber\\[-8pt]\\[-8pt]
&&\qquad\hspace*{26pt}{}\times{\frac{E[{\exp({{{\gamma(1+\varepsilon)^2}/{2}} \|{Z}\|
_{2q,R}^2})}]^{{1}/({1+\varepsilon})}}{P [{\|{Z}\|_{2q,R} \geq\sqrt
{2b_T \varepsilon} + |s| R^{d/2q}}]}}
\nonumber\\
&&\qquad\hspace*{26pt}{}\times
\exp\biggl({{{\frac{\gamma}{2}} {\frac{1+\varepsilon}{\varepsilon}} s^2
R^{d/q}}}\biggr).\nonumber
\end{eqnarray}
Choose $s= \sqrt{2b_T} \varepsilon^{3/2} R^{-d/2q}$ in (\ref{ubexp2.eq}).
$\forall\theta> 0$, $\forall\gamma> \theta$, $\forall\varepsilon>0$,
$\forall R > 0$, $\forall T > 0$,
%
\begin{eqnarray}\quad
&&E[{\exp(\theta I_{R,\tau}^{1/q})}]
\nonumber\\
&&\qquad\leq1 + C(\varepsilon) {\frac{\theta}{\gamma- \theta}}
\biggl({1+ {\frac{\sqrt{T} R^{d/2q}}{\sqrt{\alpha} b_T \varepsilon^{3/2} }}}\biggr)
\\
&&\qquad\hspace*{26pt}{}\times
{\frac{E[{\exp({{{\gamma(1+\varepsilon)^2}/{2}} \|{Z}\|
_{2q,R}^2})}]^{{1}/({1+\varepsilon})}}{P [{\|{Z}\|_{2q,R} \geq\sqrt
{2b_T \varepsilon}(1+\varepsilon)}]}}
\exp\bigl({\gamma\varepsilon^2 (1+\varepsilon) b_T}\bigr).\nonumber
\end{eqnarray}
(\ref{ubl4norm.eq}) is now obtained by the change of variable $\gamma
\leadsto\gamma/(1+\varepsilon)^2$.
\end{pf}

\textit{Step} 4.
Large deviations for $\|{Z}\|_{2q,R}$.
\begin{lemma}
\label{expnablaZ.lem}
For any $\alpha>0$ and $R>0$, let $\tau$ and $(Z_x, x\in\mathbb
{T}_R)$ be defined
as in Lemma \ref{dynkin.lem}. Let $\rho_1(\alpha,R,T):= \inf
\{\sum_{x,y\in\mathbb{T}_R} f_x G_{R,\lambda}^{-1}(x,y) f_y$;
$f$ such that $\sum_{x \in\mathbb{T}_R} f_x^{2q}=1\}$.

\begin{enumerate}
\item$\forall\alpha>0$, $\forall R >0$,
$\forall T > 0$, $\alpha\frac{b_T}{T} \leq
\rho_1(\alpha,R,T) \leq2d + \alpha\frac{b_T}{T}$.
\item$\forall\alpha>0$, $\forall\varepsilon>0$, $\forall R >0$,
$\forall T > 0$,
%
\begin{equation}
\label{LBZ.lem}
P \bigl[{\|{Z}\|_{2q,R} \geq\sqrt{b_T \varepsilon}}\bigr]
\geq\frac{1 - {1}/({b_T \varepsilon\rho_1(\alpha,R,T)})}
{ \sqrt{2 \pi b_T \varepsilon\rho_1(\alpha,R,T)}}
\exp\biggl({-\frac{b_T
\varepsilon\rho_1(\alpha,R,T)}{2}}\biggr).\hspace*{-15pt}
\end{equation}
\item$\exists C(q)$ such that
$\forall\alpha>0$, $\forall R > 0$,
$\forall T > 0$, $\forall\gamma< \rho_1(\alpha,R,T)$,
$\forall\varepsilon> 0$ such that $\gamma(1+ \varepsilon) < \rho
_1(\alpha,R,T)$,
%
\begin{eqnarray}
\label{UBZ.lem}
E\biggl[{\exp\biggl({{\frac{\gamma}{2}} \|{Z}\|_{2q,R}^2}\biggr)}\biggr]
&\leq&
{\frac{2}{\sqrt{1-{{\gamma(1+\varepsilon)}/({\rho_1(\alpha,R,T)})}}}}
\nonumber\\[-8pt]\\[-8pt]
&&\times{}
\exp\biggl({C(q) \gamma\frac{1+\varepsilon}{\varepsilon} R^{d/q} G_{R,\lambda
}(0,0)}\biggr).\nonumber
\end{eqnarray}
\end{enumerate}
\end{lemma}
\begin{pf}
1. Since $G_{R,\lambda}=(\lambda\operatorname{Id}- \triangle)^{-1}$,
\[
\rho_1(\alpha,R,T) =\inf\{{\lambda\|{f}\|^2_{2,R} - (f, \triangle
f) ; f \mbox{ such that } \|{f}\|_{2q,R}=1}\}.
\]

Taking $f=\delta_0$, we obtain that $\rho_1(\alpha,R,T)
\leq\lambda+ 2d = \alpha\frac{b_T}{T} + 2d $.
For the lower bound, note that if $\|{f}\|_{2q,R}=1$, for all $x \in
\mathbb{T}_R$, $|f_x| \leq1$, so that $\|{f}\|^{2}_{2,R} \geq\sum
_{x \in
\mathbb{T}_R} f_x^{2q} =1$. Therefore, $\rho_1(\alpha,R,T) \geq
\lambda$.


\begin{enumerate}[2.]
\item[2.] For all $(f_x, x \in\mathbb{T}_R)$, such that $\sum_x
|f_x|^{{2q}/({2q-1})}=1$,
\[
P\bigl[{\|{Z}\|_{2q,R} \geq\sqrt{b_T \varepsilon}}\bigr]
\geq P \biggl[{\sum_{x \in\mathbb{T}_R} f_x Z_x \geq\sqrt{ b_T \varepsilon
} }\biggr].
\]
$\sum_{x \in\mathbb{T}_R} f_x Z_x$ is a real centered Gaussian
variable, with
variance
\[
\sigma^2_{\alpha,R,T}(f)=
\sum_{x,y \in\mathbb{T}_R} G_{R,\lambda}(x,y) f_x f_y.
\]
Therefore, for all $(f_x, x \in\mathbb{T}_R)$, such that
$\sum_x |f_x|^{{2q}/({2q-1})}=1$,
\begin{eqnarray*}
P\bigl[{\|{Z}\|_{2q,R} \geq\sqrt{b_T \varepsilon}}\bigr]
& \geq&
\frac{\sigma_{\alpha,R,T}(f)}{\sqrt{2\pi} \sqrt{b_T \varepsilon}}
\biggl({1 - \frac{\sigma^2_{\alpha,R,T}(f)}{b_T \varepsilon}}\biggr)
\exp\biggl({- \frac{b_T \varepsilon}{2 \sigma^2_{\alpha,R,T}(f)}}\biggr),\\
& \geq&
\frac{\sigma_{\alpha,R,T}(f)}{\sqrt{2\pi} \sqrt{b_T \varepsilon}}
\biggl({1 - \frac{\rho_2(\alpha,R,T)}{b_T \varepsilon}}\biggr)
\exp\biggl({- \frac{b_T \varepsilon}{2 \sigma^2_{\alpha,R,T}(f)}}\biggr),
\end{eqnarray*}
where $\rho_2(\alpha,R,T):= \sup\{\sigma^2_{\alpha,R,T}(f)$; $f$
such that ${\sum_{x \in\mathbb{T}_R}} |f_x|^{{2q}/({2q-1})}
=1\}$.
Take the supremum over $f$, to obtain $\forall\alpha>0$, $\forall R >0$,
$\forall T >0$,
%
\begin{eqnarray}
\label{BIZ4.eq}
P\bigl[{\|{Z}\|_{2q,R} \geq\sqrt{b_T \varepsilon}}\bigr]
&\geq&
\frac{\sqrt{\rho_2(\alpha,R, T)}}{\sqrt{2 \pi b_T \varepsilon}}
\biggl({1-\frac{\rho_2(\alpha,R,T)}{b_T \varepsilon}}\biggr)
\nonumber\\[-8pt]\\[-8pt]
&&{}\times
\exp\biggl({- \frac{b_T \varepsilon}{2 \rho_2(\alpha,R,
T)}}\biggr).\nonumber
\end{eqnarray}
We are now going to prove that
$\forall\alpha>0$, $\forall R > 0$, $\forall
T > 0$,
%
\begin{equation}
\label{r1=inv_r2}
\rho_2(\alpha,R,T) = \frac{1}{\rho_1(\alpha,R,T)}.
\end{equation}
Indeed,
\begin{eqnarray*}
(G_{R,\lambda}h,h) &=& (G_{R,\lambda}h, G_{R,\lambda}^{-1} G_{R,\lambda}h)
\geq\rho_1(\alpha,R,T) \|{G_{R,\lambda}h}\|_{2q,R}^2
\\
&\geq&\rho_1(\alpha,R,T) \frac{(G_{R,\lambda}h,h)^2}
{\|{h}\|^2_{{2q}/({2q-1}),R}},
\end{eqnarray*}
where the first inequality follows from the definition of $\rho
_1(\alpha,R,T)$,
and the second one from H\"{o}lder's inequality. Therefore, for all $h$,
$(G_{R,\lambda}h,h)
\leq\frac{1}{\rho_1(\alpha,R,T)}{\|{h}\|^2_{{2q}/({2q-1}),R}}$. Taking
the supremum over $h$ yields
$ \rho_2(\alpha,R,T) \leq\frac{1}{\rho_1(\alpha,R,T)}$.
For the opposite inequality, take $f_0$ achieving the infimum in the
definition of $\rho_1(\alpha,R,T)$. Applying the Lagrange multipliers method,
it\break is easy to see that $f_0$ satisfies the equation
$G_{R,\lambda}^{-1}f_0=\rho_1(\alpha,R,T) f_0^{2q-1}$.\break Hence,
$\|{G_{R,\lambda}^{-1}f_0}\|_{{2q}/({2q-1}),R} = \rho_1(\alpha,R,T)
\|{f_0^{2q-1}}\|_{{2q}/({2q-1}),R}= \rho_1(\alpha, R,\break T)
\|{f_0}\|_{2q,R}^{2q-1}=
\rho_1(\alpha,R,T)$.
Moreover, $(G_{R,\lambda}^{-1}f_0,f_0)= \rho_1(\alpha,R,T)$ and
\[
\rho_2(\alpha,R,T) \geq
\frac{(G_{R,\lambda}^{-1}f_0,G_{R,\lambda}G_{R,\lambda}^{-1}f_0)}
{\|{G_{R,\lambda}^{-1}f_0}\|_{{2q}/({2q-1}),R}^2}
\geq\frac{\rho_1(\alpha,R,T)}{\rho_1(\alpha,R,T)^2}
= \frac{1}{\rho_1(\alpha,R,T)},
\]
which ends the proof of (\ref{r1=inv_r2}) and of (\ref{LBZ.lem}).

\item[3.] Let $M_{R,T}$ denote the median of $\|{Z}\|_{2q,R}$.
For $\gamma< \rho_1(\alpha,R,T)$, and $\varepsilon> 0$ such
that $\gamma(1+\varepsilon)< \rho_1(\alpha,R,T)$,
\begin{eqnarray*}
E\biggl[{\exp\biggl({\frac{\gamma}{2} \|{Z}\|^2_{2q,R}}\biggr)}\biggr]
&\leq& E \biggl[{\exp\biggl({\frac{\gamma(1+\varepsilon)}{2} ({\|{Z}\|_{2q,R} -
M_{R,T}})^2}\biggr)}\biggr]\\
&&{}\times\exp\biggl({\frac{\gamma}{2}\frac{1+\varepsilon}{\varepsilon}
M_{R,T}^2}\biggr).
\end{eqnarray*}
But $M_{R,T}=\mathrm{median}((\sum_x Z_x^{2q})^{1/2q})=(\mathrm{
median}(\sum
_x Z_x^{2q}))^{1/2q}$.
Moreover, it is easy to see that for any positive r.v. $X$, $\mathrm{median}(X)
\leq
2 E(X)$. Hence, using the fact that $Z_x$ is a centered Gaussian variable
with variance $G_{R,\lambda}(0,0)$,
\[
M_{R,T}^2
\leq2^{1/q} E\biggl[{\sum_{x \in\mathbb{T}_R} Z_x^{2q}}\biggr]^{1/q}
= 2^{1/q} R^{d/q} G_{R,\lambda}(0,0) E(V^{2q})^{1/q},
\]
where $V \sim\mathcal{N}(0,1)$.

On the other hand,
\begin{eqnarray*}
&&E\biggl[{\exp\biggl({{\frac{\gamma(1+\varepsilon)}{2}} ({\|{Z}\|_{2q,R} - M_{R,T}})^2}\biggr)}\biggr]
\\
&&\qquad= 1 + {\int_{0}^{\infty}} {\frac{\gamma(1+\varepsilon)}{2}}
e^{{\gamma(1+\varepsilon)u}/{2} }
P\bigl[{\bigl|{\|{Z}\|_{2q,R} - M_{R,T}}\bigr| \geq\sqrt{u}}\bigr] \,du.
\end{eqnarray*}
We now use the concentration inequalities for norms of Gaussian
processes (see, for instance, Lemma 3.1 in \cite{ledoux-talagrand}):
$\forall u > 0$,
\[
P\bigl[{\bigl|{\|{Z}\|_{2q,R} - M_{R,T}}\bigr| \geq\sqrt{u}}\bigr]
\leq2 P\bigl( V \geq\sqrt{\rho_1(\alpha, R, T) u}\bigr).
\]
Therefore, since $\gamma(1+\varepsilon) < \rho_1(\alpha, R, T)$,
\begin{eqnarray*}\hspace*{67.5pt}
&&E\biggl[{\exp\biggl({{\frac{\gamma(1+\varepsilon)}{2}} ({\|{Z}\|_{2q,R} - M_{R,T}})^2}\biggr)}\biggr]
\\
&&\qquad \leq
- 1+ 2 E\biggl[{\exp\biggl({\frac{\gamma(1+\varepsilon)}{2 \rho_1(\alpha, R,T)} V^2}\biggr)}\biggr]
\\
&&\qquad = -1 + {\frac{2}{\sqrt{1-{{\gamma(1+\varepsilon)}/({\rho
_1(\alpha,R,T)})}}}}.\hspace*{66pt}\qed
\end{eqnarray*}
\end{enumerate}
\noqed\end{pf}

\textit{Step} 5.
An upper bound for exponential moments of $I_T$ and $Q_T$.
\begin{lemma}
\label{UB.lem}
Assume that
$\log(T) \ll b_T \leq T$, and that $R$ depends on $T$ in such a way that
$\forall\alpha>0$, $b_T \gg R^{d/q} G_{R,\lambda}(0,0)$.
For all $\alpha> 0$, set
\begin{eqnarray*}
\rho_1(\alpha) &=& \liminf_{T \rightarrow\infty} \rho_1(\alpha,R,T)
\\
&=&\liminf_{T \rightarrow\infty}
\inf\biggl\{{\alpha\frac{b_T}{T} \|{f}\|^2_{2,R} + \|{\nabla f}\|
^2_{2,R} ; f \mbox{ such that } \|{f}\|_{2q,R} =1 }\biggr\}
\\
\rho_1 &=& \limsup_{\alpha\rightarrow0} \rho_1(\alpha).
\end{eqnarray*}

\begin{enumerate}
\item For any $\theta\in[0,\rho_1[$,
$\limsup_{T \rightarrow\infty} \frac{1}{b_T} \log
E [{\exp(\theta I_T^{1/q})}]
= 0$.
\item
For any $\theta\in[0,q\rho_1[$,
\[
\limsup_{T \rightarrow\infty} \frac{1}{b_T} \log
E [{\exp(\theta Q_T^{1/q}) }]
= 0.
\]
\end{enumerate}
\end{lemma}
\begin{pf}
Point 2 is a straightforward consequence of 1, since
\[
Q_T^{1/q} = \Biggl({\sum_x \prod_{i=1}^q l^{(i)}_{T}(x)}\Biggr)^{1/q}
\leq\Biggl({\prod_{i=1}^q \bigl\|{l^{(i)}_{T}}\bigr\|_{q}}\Biggr)^{1/q}
\leq\frac{1}{q} \sum_{i=1}^q \bigl\|{l^{(i)}_{T}}\bigr\|_{q},
\]
where the last inequality comes from the concavity of the $\log$ function.
Hence,
\[
E[{\exp(\theta Q_T^{1/q})}] \leq
E\biggl[{\exp\biggl({\frac{\theta}{q} \|{l_{T}}\|_{q}}\biggr)}\biggr]^q
= E\biggl[{\exp\biggl({\frac{\theta}{q} I_T^{1/q}}\biggr)}\biggr]^q.
\]

We thus focus on step 1 of Lemma \ref{UB.lem}. Let $\alpha> 0$,
and $\theta< \rho_1(\alpha)$ be fixed.
Take $\gamma$ such that $\theta< \gamma< \rho_1(\alpha)$.
Take then $\varepsilon\in\,]0;\min(\sqrt{\frac{\gamma}{\theta}-1},1)[$
such that
\[
\theta< \gamma< \gamma(1+ 2 \varepsilon) < \rho_1(\alpha).
\]
For $T$ sufficiently large ($ T \geq T_0$), $\rho_1(\alpha, R, T)
\geq\gamma(1+ 2 \varepsilon)$. Lemmas
\ref{redtorus.lem} and \ref{ubl4norm.lem}
lead to
%
\begin{eqnarray}
e^{-\alpha b_T} E[{e^{\theta I_T^{1/q}}}]
&\leq& 1 + C(\varepsilon) {\frac{\theta}{\gamma- \theta(1+\varepsilon)^2}}
\biggl({1+{\frac{\sqrt{T}R^{d/2q}}{\sqrt{\alpha} b_T}}}\biggr)
\nonumber\\[-8pt]\\[-8pt]
&&\hspace*{15.7pt}{}\times
{\frac{E [{\exp({{{\gamma}/{2}} \|{Z}\|_{2q,R}^2})}]^{
{1}/({1+\varepsilon})}}{P [{\|{Z}\|_{2q,R} \geq\sqrt{8 b_T \varepsilon}}]}}
\exp({\gamma\varepsilon^2 b_T}).\nonumber
\end{eqnarray}
By Lemma \ref{expnablaZ.lem}, for $b_T \leq T$, and $T \geq T_0$,
$\rho_1(\alpha, R, T) \geq\gamma(1+ 2 \varepsilon)$, and
\begin{eqnarray*}
&&P \bigl[{\|{Z}\|_{2q,R} \geq\sqrt{8 b_T \varepsilon}}\bigr]
\\
&&\qquad \geq \frac{1}{\sqrt{16 \pi b_T \varepsilon (2d+\alpha)}}
\biggl({1-\frac{1}{ 8 b_T \varepsilon\rho_1(\alpha,R,T)}}\biggr)
\exp\bigl(- 4 b_T \varepsilon(2d+ \alpha)\bigr), \\
&&\qquad \geq \frac{1}{\sqrt{16 \pi b_T \varepsilon (2d+\alpha)}}
\biggl({1-\frac{1}{ 8 b_T \varepsilon\gamma(1+2\varepsilon)}}\biggr)
\exp\bigl(- 4 b_T \varepsilon(2d+ \alpha)\bigr).
\end{eqnarray*}

Moreover, for $T \geq T_0$,
(\ref{UBZ.lem}) of Lemma \ref{expnablaZ.lem} yields
\begin{eqnarray*}
&&E \biggl[{\exp\biggl({\frac{\gamma}{2} \|{Z}\|_{2q,R}^2}\biggr)}\biggr]^{
{1}/({1+\varepsilon})}
\\
&&\qquad\leq
\Biggl({2 \sqrt{{\frac{1+2\varepsilon}{\varepsilon}}}}\Biggr)^{{1}/({1+\varepsilon})}
\exp\biggl({C(q) {\frac{\gamma}{\varepsilon}} R^{d/q} G_{R,\lambda}(0,0)}\biggr).
\end{eqnarray*}
Therefore, for $R^{d/q} G_{R,\lambda}(0,0) \ll b_T$, and $b_T \gg\log(T)$,
\[
\limsup_{T \rightarrow\infty} \frac{1}{b_T} \log
E[{\exp(\theta I_T^{1/q})}]
\leq\alpha+ 4 \varepsilon(2d+\alpha) + \gamma\varepsilon^2.
\]

Sending $\varepsilon$ to 0, we thus obtain that
$\forall\alpha> 0$, $\forall\theta< \rho_1(\alpha)$,
%
\begin{equation}
\label{UBIap.eq}
\limsup_{T \rightarrow\infty}
\frac{1}{b_T} \log E [{\exp(\theta I_T^{1/q})}]
\leq\alpha.
\end{equation}
Take now $\theta< \rho_1 = \limsup_{\alpha\rightarrow0} \rho
_1(\alpha)$.
Let $(\alpha_n)$ be a sequence converging to 0, such that
$\lim_{n \rightarrow\infty} \rho_1(\alpha_n) = \rho_1$. For sufficiently
large $n$, $\rho_1(\alpha_n) > \theta$, and by (\ref{UBIap.eq}),
\[
\limsup_{T \rightarrow\infty}
\frac{1}{b_T} \log E [{\exp(\theta I_T^{1/q})}] \leq\alpha_n.
\]
Point 1 is now proved by letting $n$ go to infinity.
\end{pf}

\textit{Step} 6.
Study of $\rho_1$ and $G_{R,\lambda}(0,0)$.

By Lemma \ref{UB.lem} and (\ref{ExpLBI}), we know that if $R$ is such that
$b_T \gg R^{d/q} G_{R,\lambda}(0,0)$, then $\rho_1 \leq\frac{1}{C_S^2(d)}$.
It could however happen that $\rho_1 =0$. It remains thus to determine the
values of $R$ for which $\rho_1 > 0$, and to study
the behavior of $G_{R,\lambda}(0,0)$.
\begin{lemma}[{[Behavior of $\rho_1(\alpha, R, T)$]}]
\label{ro1=cs}
Let $d \geq3$, and $q=\frac{d}{d-2}$. Let $\rho_1$ be defined as in
Lemma \ref{UB.lem}.
\begin{enumerate}
\item
Assume that R depends on $T$ in such a way that $\forall\alpha>0$,
$\lambda R^2 \gg1$.
Then $\rho_1 \geq\frac{1}{C_s^2(d)}$.
\item
Assume that R depends on $T$ in such a way that $\lim_{T \rightarrow
\infty}
\lambda R^2 = l(\alpha) \in\,]0;\break+\infty[$.
Then there exists a constant $C$ such that $\forall\alpha> 0$, $\rho
_1(\alpha) >
C \min(1,l(\alpha))$.
\end{enumerate}
\end{lemma}
\begin{pf}
Let $f_0 \in l_{2q}(\mathbb{T}_R)$ achieve
the minimum in the definition of $\rho_1(\alpha,\break R,T)$. $f_0$ is
viewed as
a periodic function on $\mathbb{Z}^d$, and by definition
\[
\rho_1(\alpha,R,T)
= \lambda\|{f_0}\|^2_{2,R} + \|{\nabla f_0}\|^2_{2,R} ;\qquad
\|{f_0}\|_{2q,R} = 1.
\]

Let $0<r<R$, and define
\[
\mathcal{C}_{r,R} = \bigcup_{i=1}^d
\{{x\in\mathbb{Z}^d ; 0 \leq x_i \leq r \mbox{ or } R-r \leq x_i
\leq R}\}
.
\]
Then one can find $a \in\mathbb{Z}^d$ such that $\sum_{x \in
\mathcal{C}_{r,R}}
f_0^{2q}(x-a)
\leq\frac{2dr}{R}$. Indeed, on one hand,
\begin{eqnarray*}
\sum_{a \in[0,R]^d} \sum_{x \in\mathcal{C}_{r,R}} f_0^{2q}(x-a)
& = & \sum_{x \in\mathcal{C}_{r,R}} \sum_{a \in[0,R]^d} f_0^{2q}(x-a)
\\
& = & \sum_{x \in\mathcal{C}_{r,R}} \sum_{x \in\mathbb{T}_R}
f_0^{2q}(x)
= \operatorname{card}(\mathcal{C}_{r,R}) \leq2d r R^{d-1}.
\end{eqnarray*}
On the other hand,
\[
\sum_{a \in[0,R]^d} \sum_{x \in\mathcal{C}_{r,R}} f_0^{2q}(x-a)
\geq R^d \inf_{a \in[0;R]^d} \sum_{x \in\mathcal{C}_{r,R}} f_0^{2q}(x-a)
.
\]

Set $f_{0,a}(x)\triangleq f_0(x-a)$. $f_{0,a}$
is a periodic function of period $R$. Note that
$\|{\nabla f_{0,a}}\|_{2,R} = \|{\nabla f_0}\|_{2,R}$,
$\|{f_{0,a}}\|_{2q,R} = \|{f_0}\|_{2q,R}$, and that
$\|{f_{0,a}}\|_{2,R} = \|{f_0}\|_{2,R}$. We can therefore assume without
loss of generality, that $f_0$ achieving the minimum in the definition
of $\rho_1(\alpha,R,T)$, satisfies also
\[
\sum_{x \in\mathcal{C}_{r,R}} f_0^{2q}(x) \leq\frac{2dr}{R}.
\]
Let $\psi\dvtx \mathbb{Z}^d \mapsto[0,1]$ a truncature function satisfying
\[
\cases{
\psi(x) = 0, &\quad if $x \notin[0;R]^d$;
\cr
\psi(x) = 1, &\quad if $x \in[0;R]^d \setminus\mathcal{C}_{r,R}$;
\cr
|\nabla_i \psi(x)| \leq\dfrac{1}{r},
&\quad $\forall x \in\mathbb{Z}^d, \forall i \in\{{1,\ldots, d}\}$.}
\]
Fix $\varepsilon> 0$, and take $r=\frac{\varepsilon R}{2d}$.
By definition, for $q=\frac{d}{d-2}$,
\[
\frac{1}{C_s^2(d)}
\leq\frac{\|{\nabla(\psi f_0)}\|_2^2}{\|{\psi f_0}\|_{2q}^2}.
\]
Regarding the denominator,
%
\begin{equation}
\label{nor4psif}
\|{\psi f_0}\|^{2q}_{2q} \geq\sum_{x \in[0;R]^d} f_0^{2q}(x) -
\sum_{x \in\mathcal{C}_{r,R}} f_0^{2q}(x) \geq1- \frac{2dr}{R}= 1-
\varepsilon.
\end{equation}
It remains to control $\|{\nabla(\psi f_0)}\|_2$,
%
\begin{eqnarray}\label{nor2gradpsif}
\|{\nabla(\psi f_0)}\|_2^2
& = & \sum_{x \in[0;R]^d} \sum_{i=1}^d \bigl({\nabla_i \psi(x)
f_0(x+e_i) + \psi(x) \nabla_i f_0(x)}\bigr)^2
\nonumber\\
& = & \sum_{x \in[0;R]^d} \sum_{i=1}^d ({\nabla_i \psi(x)})^2
f_0^2(x+e_i) + \psi^2(x) ({\nabla_i f_0(x)})^2
\nonumber\\
& &{} + 2 \sum_{x \in[0;R]^d} \sum_{i=1}^d
\nabla_i \psi(x) \psi(x) f_0(x+e_i) \nabla_i f_0(x)
\nonumber\\[-8pt]\\[-8pt]
& \leq& \frac{d}{r^2} \|{f_0}\|^2_{2,R} + \|{\nabla f_0}\|_{2,R}^2
+ \frac{2 \sqrt{d}}{r} \|{f_0}\|_{2,R} \|{\nabla f_0}\|_{2,R}
\nonumber\\
& \leq& \|{\nabla f_0}\|_{2,R}^2 (1+\varepsilon)
+ \frac{d}{r^2} \|{f_0}\|^2_{2,R} (1 + 1/\varepsilon).
\nonumber\\
& \leq& (1+\varepsilon) \max\biggl(1, \frac{d}{\lambda r^2 \varepsilon}\biggr) \rho
_1(\alpha, R, T).\nonumber
\end{eqnarray}
It follows from (\ref{nor4psif}) and (\ref{nor2gradpsif})
that $\forall\varepsilon
\in\,]0;1[$, $\forall\alpha>0$, $\forall T >0$,
%
\begin{equation}
\label{rho1.eq}
\frac{1}{C_S^2(d)} \leq\frac{1 + \varepsilon}{(1-\varepsilon)^{1/q}}
\max\biggl(1, \frac{4 d^3}{\varepsilon^3} \frac{1}{\lambda R^2}\biggr) \rho
_1(\alpha,
R, T).
\end{equation}

\textit{Case} 1.
Since $R$ is such that $b_T \gg
\frac{T}{R^2}$, $\forall\varepsilon> 0$, $\forall\alpha> 0$, $\rho
_1(\alpha)
\geq\frac{1}{C_S^2(d)} \frac{(1-\varepsilon)^{1/q}}{1+\varepsilon}$.
Hence, letting
$\varepsilon$ go to $0$, $\forall\alpha> 0$, $\rho_1(\alpha) \geq
\frac
{1}{C_S^2(d)}$,
so that $\rho_1 \geq\frac{1}{C_S^2(d)}$.

\textit{Case} 2.
Take $\varepsilon= 1/2$ in (\ref{rho1.eq}), and let $l(\alpha) = \lim_{T
\rightarrow
\infty} \lambda R^2$. Then $\forall\alpha> 0$,
\[
\rho_1(\alpha)
\geq\frac{2^{1-1/q}}{3} \frac{1}{C_s^2(d)} \min\biggl(1, \frac{l(\alpha
)}{32 d^3}\biggr)
\geq C \min(l(\alpha),1).
\]
\upqed\end{pf}
\begin{lemma}[{[Behavior of $G_{R,\lambda}(0,0)$]}]
\label{Green}
Assume that $d \geq3$, that $\lambda\ll1$, and that
$R$ depends on $T$ in such a way that $\lambda R^d \gg1$. Then
$\lim_{T \rightarrow\infty} G_{R,\lambda}(0,0) = G_d(0,0)$, where
$G_d(0,0)$ is the expected amount of time the simple random walk on
$\mathbb{Z}^d$
spends on site $0$.
\end{lemma}
\begin{pf}
Let $p^R_t(x,y)$ be the transition probability of $X^{(R)}_t$. Then
\[
G_{R,\lambda}(0,0) = \int_0^{\infty} \exp(-\lambda t) p^R_t(0,0)\,
dt.
\]
It follows from Nash inequality (see, for instance,
Theorems 2.3.1 and 3.3.15 in \cite{saloff}) that
there exists a constant $C(d)$ such that $\forall R>0$, $\forall t > 0$,
\[
\biggl|{p^R_t(0,0) - \frac{1}{R^d}}\biggr|
\leq\frac{C(d)}{t^{d/2}}.
\]
Therefore, $\forall S >0$,
\begin{eqnarray*}
&&\int_S^{+\infty} \exp(-\lambda t) p^R_t(0,0) \,dt
\\
&&\qquad\leq
\frac{1}{R^d} \int_0^{\infty} \exp(-\lambda t) \,dt
+ \int_S^{+\infty} \frac{C(d)}{t^{d/2}} \,dt\\
&&\qquad\leq\frac{1}{\lambda R^d} + \frac{C(d)}{S^{{d}/{2}-1}}.
\end{eqnarray*}
Thus, when $\lambda R^d \gg1$, and $S \gg1$,
%
\begin{equation}
\label{Ggrandt}
\lim_{T \rightarrow\infty} \int_S^{+\infty} \exp(-\lambda t)
p^R_t(0,0) \,dt = 0.
\end{equation}
For the values of $t$ less than $S$,
\begin{eqnarray*}
p^R_t(0,0)& = & P_0\bigl(X^{(R)}_t=0\bigr)
\\
& \leq& P_0 \biggl[{X^{(R)}_t=0; \sup_{s \leq S} \|{X_s}\| \leq\frac{R}{2}}\biggr]
+ P_0 \biggl[{\sup_{s \leq S} \|{X_s}\| \geq\frac{R}{2}}\biggr]
\\
& = & P_0 \biggl[{X_t=0; \sup_{s \leq S} \|{X_s}\| \leq\frac{R}{2}}\biggr]
+ P_0 \biggl[{\sup_{s \leq S} \|{X_s}\| \geq\frac{R}{2}}\biggr]
\\
& \leq& P_0 [{X_t=0}] + C(d) \exp\biggl(-\frac{R^2}{C(d)S}\biggr).
\end{eqnarray*}
The third equality comes from the fact that as long as $X$ does not
exit a
ball of radius $R/2$, then $X$ and $X^{(R)}$ are the same. The fourth
one follows
from standard results on simple random walks.
Thus,
\[
\int_0^{S} \exp(-\lambda t) p^R_t(0,0) \,dt
\leq\int_0^{\infty} p_t(0,0) \,dt + C(d) S \exp\biggl(-\frac{R^2}{C(d)S}\biggr).
\]
On the other hand, $p^R_t(0,0) = P_0(X^{(R)}_t=0) \geq p_t(0,0)$, so that
\begin{eqnarray*}
\int_0^{S} \exp(-\lambda t) p^R_t(0,0) \,dt
& \geq& \int_0^{S} p_t(0,0) \,dt - \int_0^{S} \bigl(1-\exp(-\lambda t)\bigr)\,
dt
\\
& = & \int_0^{S} p_t(0,0) \,dt + \frac{\exp(-\lambda S)-1+\lambda
S}{\lambda}.
\end{eqnarray*}
Hence, if $S$ is chosen so that $S \gg1$,
$S \ll R^2 / (\log(R))^{1+\varepsilon}$,
and $\lambda S^2 \ll1$,
%
\begin{equation}
\label{Gpetitt}
\lim_{T \rightarrow\infty} \int_0^{S} \exp(-\lambda t) p^R_t(0,0)
\,dt
= \int_0^{\infty} p_t(0,0) \,dt = G_d(0,0).
\end{equation}
Now, for $\lambda\ll1$, and $\lambda R^d \gg1$ (which implies $R \gg
1$), one
can always choose $S$ such that $1 \ll S \ll\min(R^2/(\log
(R))^{1+\varepsilon},
1/\sqrt{\lambda})$.
For such a choice of $S$, it follows from (\ref{Ggrandt}) and (\ref
{Gpetitt}) that
\[
\lim_{T \rightarrow\infty} G_{R,\lambda}(0,0)
= G_d(0,0) < \infty \qquad\mbox{for } d \geq3.
\]
\upqed\end{pf}

\textit{Step} 7.
End of proof of Proposition \ref{SILT4}.

Choose $R$ such that
\[
\frac{T}{R^2} \ll b_T,\qquad b_T \gg R^{d/q}.
\]
Then, on one hand, $\forall\alpha> 0$, $\lambda b_T \ll R^2$, and
$\rho_1 \geq\frac{1}{C_S^2(d)}$ by 1. of Lemma \ref{ro1=cs}. On the other
hand, $\lambda R^d = \alpha\frac{b_T}{T} R^d \gg\alpha\frac
{b_T}{T} R^2
\gg1$. Hence, by Lemma \ref{Green}, $G_{R,\lambda}(0,0) \simeq G_d(0,0)$
and it follows from Lemma \ref{UB.lem} that $\rho_1 \leq\frac{1}{C_S^2(d)}$.
Therefore, for such a choice of $R$, $\rho_1= \frac{1}{C_s^2(d)}$ and
\begin{eqnarray*}
\forall\theta&\in&\biggl[0;\frac{1}{C_s^2(d)}\biggr[\qquad
\liminf_{T \rightarrow\infty} \frac{1}{b_T} \log E [{\exp(\theta
I_T^{1/q})}] = 0,
\\
\forall\theta&\in&\biggl[0;\frac{q}{C_s^2(d)}\biggr[
\qquad\liminf_{T \rightarrow\infty}
\frac{1}{b_T} \log E [{\exp(\theta Q_T^{1/q})}] = 0.
\end{eqnarray*}

The best choice for $R$ corresponds to $T/R^2=R^{d/q}=R^{d-2}$,
i.e., $R^d=T$, leading
to $b_T \gg T^{1-2/d}=T^{1/q}$.


\section{Large and very large deviations lower bounds}
\label{LBGD}
The aim of this section is to prove the lower bounds in Theorems \ref
{ILT4VLD} and
\ref{SILT4LD}. We have actually the following result.
\begin{prop}\label{SILT4LB}
1. Lower bound for $I_T$.

Assume that $d \geq3$, $q=\frac{d}{d-2}$, and $T \gg b_T \gg1$.
%
\begin{equation}
\label{SILT4LB.eq}
\liminf_{T \rightarrow\infty} \frac{1}{b_T} \log P [{I_T \geq b_T^q}]
\geq- \frac{1}{C_S^2(d)}.
\end{equation}

\begin{enumerate}[2.]
\item[2.] Lower bound for $Q_T$.

Assume that $d=4$ and $q=2$, or $d=3$ and $q=3$, and that $1\ll b_T
\ll T$.
%
\begin{equation}
\label{MILT4LB.eq}
\liminf_{T \rightarrow\infty} \frac{1}{b_T} \log P [{Q_T \geq b_T^q}]
\geq- \frac{q}{C_S^2(d)}.
\end{equation}
\end{enumerate}
\end{prop}
\begin{pf*}{Proof of (\protect\ref{SILT4LB.eq})}
Fix $M > 0$.
Let $T_0$ be such that for all $T \geq T_0$, $\frac{T}{b_T} > M$. For
$T \geq T_0$,
\[
P[{I_T \geq b_T^q }]
\geq P[{I_{Mb_T} \geq b_T^q }]
\geq P\biggl[{\biggl\|{\frac{l_{Mb_T}}{Mb_T}}\biggr\|_q
\geq\frac{1}{M}}\biggr].
\]
The function $\mu\in\mathcal{F}\mapsto\|{\mu}\|_q
=\sup_{f; \|{f}\|_{q'}=1} \sum_x \mu(x) f(x)$ is lower semicontinuous
in $\tau$-topology, so that
$\forall t >0$, $\{{\mu\in\mathcal{F}, \|{\mu}\|_q >t}\}$ is an open
subset of $\mathcal{F}$. Therefore, $\forall\varepsilon>0$,
\begin{eqnarray*}
&&\liminf_{T \rightarrow\infty} \frac{1}{Mb_T} \log
P \biggl[{\biggl\|{\frac{l_{Mb_T}}{Mb_T}}\biggr\|_q
\geq\frac{1}{M}}\biggr]\\
&&\qquad \geq\liminf_{T \rightarrow\infty}
\frac{1}{Mb_T} \log
P \biggl[{\biggl\|{\frac{l_{Mb_T}}{Mb_T}}\biggr\|_q > \frac{1-\varepsilon}{M}}\biggr]
\\
&&\qquad \geq - \inf\biggl\{{\|{\nabla f}\|^2_2 ; \|{f}\|_2=1, \|{f}\|
_{2q}^2 > \frac{1-\varepsilon}{M}}\biggr\}.
\end{eqnarray*}

We have thus proved that $\forall M >0$,
$\forall\varepsilon>0$,
\[
\liminf_{T \rightarrow\infty} \frac{1}{b_T} \log
P [{I_T \geq b_T^q}] \geq- M \rho_3\biggl({\frac{1-\varepsilon}{M}}\biggr),
\]
where $\rho_3(y) := \inf\{{\|{\nabla f}\|^2_2 ; \|{f}\|
^2_{2q} > y, \|{f}\|_2=1}\}$.
To end the proof of (\ref{SILT4LB.eq}), it remains to show that when
$q=\frac{d}{d-2}$, $\forall
y > 0$,
%
\begin{equation}
\label{ro3=Sob}
\inf_{M > 0} M \rho_3(y/M) = \frac{y}{C_S^2(d)}.
\end{equation}
But, if $q=\frac{d}{d-2}$,
$\forall y > 0$,
%
\begin{eqnarray}
\inf_{M > 0} M \rho_3(y/M)
& = & y \inf_{M > 0} M \rho_3(1/M)
\\
& = & y \inf_{M > 0} \inf_{f} \biggl\{{M \|{\nabla f}\|_2^2 ; \|{f}\|
_2=1, \|{f}\|_{2q}^2 > \frac{1}{M} }\biggr\}
\\
& = & y \inf_{f;\|{f}\|_2=1} \inf_{M > 0} \biggl\{{M \|{\nabla f}\|_2^2; M
> \frac{1}{ \|{f}\|_{2q}^2}}\biggr\}
\\
\label{ro3.eq}
& = & y \inf_{f;\|{f}\|_2=1} \biggl\{{\frac{\|{\nabla f}\|_2^2}{\|{f}\|
_{2q}^2}}\biggr\} ;
\\
& = & {\frac{y}{C_S^2(d)}}.
\end{eqnarray}
\upqed\end{pf*}
\begin{pf*}{Proof of (\protect\ref{MILT4LB.eq})}
The proof of (\ref{MILT4LB.eq}) cannot be done as the proof of
(\ref{SILT4LB.eq}), since the function $(\mu_1,\ldots,\mu_q)
\mapsto\sum_{x \in\mathbb{Z}^d} \mu_1(x) \cdots\mu_q(x)$ is not lower
semicontinuous in the product of $\tau$-topology.

Let $\varepsilon> 0$ be fixed. Let $h$ be a
function approaching the
infimum in the definition of $C_S(d)$, i.e., $h$ is such that
\[
\|{\nabla h}\|_2^2 \leq\frac{\|{h}\|_{2q}^2}{C_S(d)^2} (1+ \varepsilon)
,\qquad q = \frac{d}{d-2}.
\]
Dividing $h$ by its $l_2$-norm if necessary, we may and we do assume that
$\|{h}\|_2=1$.

Set $\eta= 2^{({q+1})/{q}} \varepsilon^{1/q}$,
and $M=\frac{1}{(2-(1+\eta)^q)^{1/q} \|{h}\|^2_{2q}}$ [$\varepsilon$ is
chosen\vspace*{-2pt} small enough in order that $M$ is strictly positive; actually,
one has to choose $\varepsilon< \varepsilon_0(q)= (2^{1/q}-1)^q 2^{-(q+1)}$].
For $T$ large enough, $T \geq Mb_T$, and
\[
P [{Q_T \geq b_T^q}] \geq P[{Q_{Mb_T} \geq b_T^q}].
\]
Assume that $\forall i \in\{{1,\ldots,q}\}$,
$\|{\frac{l^{(i)}_{Mb_T}}{Mb_T} - h^2}\|_{q} < \eta\|{h}\|_{2q}^2$.
Then
\begin{eqnarray*}
&&\biggl|{\frac{Q_{Mb_T}}{(Mb_T)^q} - \|{h}\|_{2q}^{2q}}\biggr|
\\
&&\qquad = \Biggl|{\sum_{x \in\mathbb{Z}^d} \prod_{1}^q \frac
{l^{(i)}_{Mb_T}(x)}{Mb_T} - h^{2q}(x)}\Biggr|
\\
&&\qquad \leq
\sum_{x \in\mathbb{Z}^d} \sum_{j=1}^q \Biggl({\prod_{i=1}^{j-1} h^2(x)}\Biggr)
\biggl|{\frac{l^{(j)}_{Mb_T}(x)}{Mb_T} - h^2(x)}\biggr|
\Biggl({\prod_{l=j+1}^q \frac{l^{(l)}_{Mb_T}(x)}{Mb_T}}\Biggr)
\\
&&\qquad \leq \sum_{j=1}^q \biggl\|{\frac{l^{(j)}_{Mb_T}}{Mb_T} - h^2}\biggr\|_q
\|{h}\|_{2q}^{2(j-1)} \prod_{l=j+1}^q \biggl\|{\frac
{l^{(l)}_{Mb_T}}{Mb_T}}\biggr\|_q
\\
&&\qquad \leq
\eta\|{h}\|_{2q}^{2q} \sum_{j=1}^q (1+\eta)^{q-j}
= \eta\|{h}\|_{2q}^{2q} \frac{(1+\eta)^q -1}{\eta}
\\
&&\qquad= [{(1+\eta)^q -1}] \|{h}\|_{2q}^{2q}.
\end{eqnarray*}

Therefore,
$Q_{Mb_T} \geq b_T^q M^q \|{h}\|_{2q}^{2q} (2-(1+\eta)^q) = b_T^q$,
by the choice of $M$.

Hence, for $T$ large enough,
%
\begin{eqnarray}
\label{LDLBQ.eq1}
P [{Q_T \geq b_T^q}]
&\geq& P\biggl[{\forall i \in\{{1,\ldots,q}\},  \biggl\|{\frac
{l^{(i)}_{Mb_T}}{Mb_T} - h^2}\biggr\|_{q} < \eta\|{h}\|_{2q}^2}\biggr]
\nonumber\\[-8pt]\\[-8pt]
&=& P\biggl[{\biggl\|{\frac{l_{Mb_T}}{Mb_T} - h^2}\biggr\|_{q} < \eta\|{h}\|_{2q}^2}\biggr]^q
.\nonumber
\end{eqnarray}
But,
\begin{eqnarray*}
\biggl\|{\frac{l_{Mb_T}}{Mb_T} - h^2}\biggr\|^q_{q}
& = & \sum_{x \in\mathbb{Z}^d} \biggl({\frac{l_{Mb_T}(x)}{Mb_T} - h^2(x)}\biggr)^q
\\
& = & \sum_{x \in\mathbb{Z}^d} \sum_{j=0}^q (-1)^{q-j}
C^q_j \frac{l^j_{Mb_T}(x)}{(Mb_T)^j} h^{2(q-j)}(x)
\\
& = & \biggl\|{\frac{l_{Mb_T}}{Mb_T}}\biggr\|_{q}^{q} + (-1)^q \|{h}\|_{2q}^{2q}
- F_q\biggl({\frac{l_{Mb_T}}{Mb_T}}\biggr),
\end{eqnarray*}
where the function $F_q$ is defined by $F_q(\mu) =
\sum_{j=1}^{q-1} (-1)^{q+1-j} C^q_j
\sum_x \mu^{j}(x)\times\break h^{2(q-j)}(x)$. Hence, for $T$ large enough,
%
\begin{eqnarray}
&&
P[{Q_T \geq b_T^q}]^{1/q}
\nonumber\\[-8pt]\\[-8pt]
&&\qquad \geq P\biggl[{F_q \biggl({\frac{l_{Mb_T}}{Mb_T}}\biggr) > \biggl\|{\frac
{l_{Mb_T}}{Mb_T}}\biggr\|_q^q +\bigl((-1)^q - \eta^q\bigr) \|{h}\|_{2q}^{2q}}\biggr]
\nonumber\\
&&\qquad \geq
P\biggl[F_q \biggl({\frac{l_{Mb_T}}{Mb_T}}\biggr) > \biggl\|{\frac{l_{Mb_T}}{Mb_T}}\biggr\|_q^q
+\bigl((-1)^q - \eta^q\bigr) \|{h}\|_{2q}^{2q} ;\nonumber\\[-8pt]\\[-8pt]
&&\hspace*{129.2pt} \biggl\|{\frac{l_{Mb_T}}{Mb_T}}\biggr\|
^q_q < \biggl(1+\frac{\eta^q}{2}\biggr)\|{h}\|_{2q}^{2q}\biggr]
\nonumber\\
\label{LDLBQ.eq2}
&&\qquad \geq
P\biggl[{F_q \biggl({\frac{l_{Mb_T}}{Mb_T}}\biggr) > \biggl(1+(-1)^q - \frac{\eta^q}{2}\biggr) \|
{h}\|_{2q}^{2q} }\biggr] \nonumber\\[-8pt]\\[-8pt]
&&\qquad\quad{}- P \biggl[{\biggl\|{\frac{l_{Mb_T}}{Mb_T}}\biggr\|^q_q \geq\biggl(1+\frac
{\eta^q}{2}\biggr)\|{h}\|_{2q}^{2q} }\biggr].\nonumber
\end{eqnarray}
The second term is controlled by the large deviations upper bound for
$I_T$, and we have
%
\begin{eqnarray}
\label{LDLBQ.eq3}
&&\limsup_{T \rightarrow\infty} {\frac{1}{b_T}}
\log P \biggl[{\biggl\|{{\frac{l_{Mb_T}}{Mb_T}}}\biggr\|^q_q \geq\biggl(1+{\frac{\eta
^q}{2}}\biggr)\|{h}\|_{2q}^{2q} }\biggr]
\nonumber\\[-8pt]\\[-8pt]
&&\qquad\leq- M {\frac{(1+ {{\eta^q}/{2}})^{1/q}}{C^2_S(d)}} \|{h}\|_{2q}^{2}
= - {\frac{(1+ {{\eta^q}/{2}})^{1/q}}{(2-(1+\eta)^q)^{1/q}
C^2_S(d)}},\nonumber
\end{eqnarray}
by the choice of $M$.

On the other hand, the function $\mu\in\mathcal{F}\mapsto F_q(\mu)$
is lower semicontinuous
in $\tau$-topology.
Indeed:
\begin{itemize}
\item For $d=4$ and $q=2$, $F_2(\mu) = 2 \sum_x \mu(x) h^2(x)$ is
continuous.
\item For $d=3$ and $q=3$, $F_3(\mu)
= 3 \sum_x \mu^2(x) h^2(x) - 3 \sum_x \mu(x) h^4(x)
= 3\times\break \sup_{g; \|{g}\|_2=1} \{{\sum_x \mu(x) h(x) g(x)}\}^2
- 3 \sum_x \mu(x) h^4(x)
$ is lower semicontinuous.
\end{itemize}
Using the large deviations lower bound in $\mathcal{F}$ for
$\frac{l_{Mb_T}}{Mb_T}$, we get that
%
\begin{eqnarray}\qquad
&&\liminf_{T \rightarrow\infty}
\frac{1}{b_T} \log
P\biggl[{F_q \biggl({\frac{l_{Mb_T}}{Mb_T}}\biggr) > \biggl(1+(-1)^q - \frac{\eta^q}{2}\biggr) \|
{h}\|_{2q}^{2q} }\biggr]
\nonumber\\[-8pt]\\[-8pt]
&&\qquad\geq
- M \inf\biggl\{{\|{\nabla g}\|_2^2 ; \|{g}\|_2 =1, F_q(g^2) >
\biggl(1+(-1)^q - \frac{\eta^q}{2}\biggr)
\|{h}\|_{2q}^{2q}}\biggr\}.\nonumber
\end{eqnarray}
Note that:
\begin{itemize}
\item For $d=4$ and $q=2$, $F_2(h^2) = 2 \|{h}\|_4^4 >
(1+(-1)^2 - \frac{\eta^2}{2}) \|{h}\|_{4}^{4}$.
\item For $d=3$ and $q=3$, $F_3(h^2) = 0
> (1+(-1)^3 - \frac{\eta^3}{2}) \|{h}\|_{6}^{6}$.
\end{itemize}
Therefore, in any case,
%
\begin{eqnarray}
\label{LDLBQ.eq4}
&&\liminf_{T \rightarrow\infty}
\frac{1}{b_T} \log
P\biggl[{F_q \biggl({{\frac{l_{Mb_T}}{Mb_T}}}\biggr) > \biggl(1+(-1)^q - {\frac{\eta^q}{2}}\biggr)
\|{h}\|_{2q}^{2q} }\biggr]
\nonumber\\
&&\qquad\geq
- M \|{\nabla h}\|_2^2
= - {\frac{\|{\nabla h}\|_2^2}{(2-(1+\eta)^q)^{{1}/{q}} \|{h}\|_{2q}^2}}
\\
&&\qquad
\geq- {\frac{1+\varepsilon}{C_S^2(d) (2-(1+\eta)^q)^{{1}/{q}}}}
,\nonumber
\end{eqnarray}
by the choice of $M$ and $h$. Putting (\ref{LDLBQ.eq2}), (\ref{LDLBQ.eq3})
and (\ref{LDLBQ.eq4}) together, we get that
%
\begin{equation}
{\frac{1}{q}} \liminf_{T \rightarrow+ \infty} {\frac{1}{b_T}}
\log P[{Q_T \geq b_T^q}]
\geq- {\frac{\min({1+\varepsilon;(1+ {\eta
^q}/{2})^{1/q}})}{(2-(1+\eta)^q)^{1/q} C^2_S(d)}}.
\end{equation}
But for $\varepsilon\in\,]0;1]$,
$(1+\varepsilon)^q = \sum_{k=0}^{q} C_q^k \varepsilon^{k}
\leq1 + \varepsilon\sum_{k=1}^{q} C_q^k = 1 + \varepsilon(2^q-1) < 1 +
\varepsilon2^q =
1+ \frac{\eta^q}{2}$. We have thus proved that
$\forall\varepsilon\in\,]0;1\wedge\varepsilon_0(q)[$,
%
\begin{equation}
\liminf_{T \rightarrow+ \infty} {\frac{1}{b_T}}
\log P[{Q_T \geq b_T^q}]
\geq- {\frac{q(1+\varepsilon)}{C_S^2(d) (2-(1+\eta)^q)^{1/q}}}.
\end{equation}
(\ref{MILT4LB.eq}) is then obtained by letting $\varepsilon$ go to zero.
\end{pf*}

\section{Large deviations upper bound}
\label{UBGD}

The only thing that remains to prove now is the upper bound in Theorem
\ref{SILT4LD}.

Let $\alpha> 0$ and $A > 0$ to be chosen later.
We take here
\[
\lambda= \alpha\frac{T^{1/q}}{T} ;\qquad R^d = A T.
\]
Let $\tau$ be an exponential time with parameter $\lambda$, independent
on the random walk. Exactly as in (\ref{UB_I.eq}), $\forall s > 0$,
$\forall\varepsilon> 0$,
%
\begin{eqnarray}
&&\exp(- \alpha T^{1/q}) P[{I_T \geq T y}]
\nonumber\\
&&\qquad\leq P[{I_{R,\tau} \geq Ty}]
\nonumber\\
&&\qquad\leq{\frac{E[{(1+ {Z_0}/{s}); \|{Z+s \mathbh{1}}\|_{2q,R} \geq
\sqrt{2} T^{1/2q} (y+\varepsilon)^{1/2q}}]}{P[{\|{Z+s \mathbh{1}}\|
_{2q,R} \geq\sqrt{2} T^{1/2q} \varepsilon^{1/2q}}]}}
\\
&&\qquad\leq
E\biggl[{\biggl(1+ {\frac{Z_0}{s}}\biggr)^{({1+\varepsilon})/{\varepsilon}}}\biggr]^{
{\varepsilon}/({1+\varepsilon})}
\nonumber\\
&&\qquad\quad{}\times
{\frac{P[{\|{Z}\|_{2q,R} \geq\sqrt{2} T^{1/2q} (y+\varepsilon)^{1/2q}
- s R^{d/2q}}]^{{1}/({1+\varepsilon})}}{P[{\|{Z}\|_{2q,R} \geq\sqrt
{2} T^{1/2q} \varepsilon^{1/2q} + s R^{d/2q}}]}}.\nonumber
\end{eqnarray}
We now choose $s R^{d/2q} = \sqrt{2} T^{1/2q} \varepsilon^{1/2q}$, i.e.,
$s = \sqrt{2} A^{-1/2q} \varepsilon^{1/2q}$.
%
\begin{eqnarray}
\label{LDUB}
&&P[{I_T \geq T y}]
\nonumber\\
&&\qquad\leq\exp( \alpha T^{1/q})
E\biggl[{\biggl(1+ {\frac{Z_0}{s}}\biggr)^{({1+\varepsilon})/{\varepsilon}}}\biggr]^{
{\varepsilon
}/({1+\varepsilon})}
\\
&&\qquad\quad{}\times
{\frac{P[{\|{Z}\|_{2q,R} \geq\sqrt{2} T^{1/2q} ((y+\varepsilon)^{1/2q}
- \varepsilon^{1/2q}) }]^{{1}/({1+\varepsilon})}}{P[{\|{Z}\|_{2q,R} \geq
2 \sqrt{2} T^{1/2q} \varepsilon^{1/2q}}]}}.\nonumber
\end{eqnarray}
Using
the fact that $Z_0$ is a centered Gaussian variable with variance
$G_{R,\lambda}(0,0)$,
we obtain that $\forall\varepsilon> 0$,
\begin{eqnarray*}
E\biggl[{\biggl(1+ \frac{Z_0}{s}\biggr)^{({1+\varepsilon})/{\varepsilon}}}\biggr]^{
{\varepsilon
}/({1+\varepsilon})}
&\leq& C(\varepsilon) \biggl({1+\frac{\sqrt{G_{R,\lambda}(0,0)}}{s}}\biggr)
\\
&\leq& C(\varepsilon) \bigl({1+ \sqrt{G_{R,\lambda}(0,0)} A^{1/2q}}\bigr).
\end{eqnarray*}
But, $\lambda R^d= \alpha A T^{1/q} \gg1$, so that $\limsup_{T
\rightarrow\infty}
G_{R,\lambda}(0,0) < \infty$ by Lemma \ref{Green}. Therefore,
$\forall
\varepsilon>0$, $\forall\alpha> 0$,
$\forall A > 0$,
\[
\limsup_{T \rightarrow\infty} \frac{1}{T^{1/q}} \log
E\biggl[{\biggl(1+ \frac{Z_0}{s}\biggr)^{({1+\varepsilon})/{\varepsilon}}}\biggr]^{
{\varepsilon
}/({1+\varepsilon})}
= 0.
\]
Let us treat the numerator of the ratio appearing in the left-hand side of
(\ref{LDUB}). Using again that
\begin{eqnarray*}
M_{R,T}&=&\mathrm{median}(\|{Z}\|_{2q,R})
\leq 2^{1/2q} E\biggl[{\sum_x Z_x^{2q}}\biggr]^{1/2q}
\\
& \leq& C(q) R^{d/2q} G_{R,\lambda}(0,0)^{1/2}
\\
&\sim& C(q) A^{1/2q} T^{1/2q} G_{d}(0,0)^{1/2},
\end{eqnarray*}
we conclude that there exists a constant $C(q)$ such that $\forall
\alpha>0$,
$\forall A >0$, for $T$ large enough, $\forall\varepsilon>0$,
%
\begin{eqnarray}\hspace*{18pt}
&&P\bigl[{\|{Z}\|_{2q,R} \geq\sqrt{2} T^{1/2q} \bigl((y+\varepsilon)^{1/2q} -
\varepsilon^{1/2q}\bigr)}\bigr]
\nonumber\\
&&\qquad\leq P\bigl[{\|{Z}\|_{2q,R} - M_{R,T} \geq\sqrt{2} T^{1/2q} \bigl((y+\varepsilon
)^{1/2q} - \varepsilon^{1/2q} -C(q) A^{1/2q} \bigr)}\bigr]
\\
&&\qquad\leq2 \exp\bigl({- T^{1/q} \rho_1(\alpha,R,T) \bigl((y+\varepsilon)^{1/2q} -
\varepsilon^{1/2q} -C(q) A^{1/2q} \bigr)_+^2}\bigr).\nonumber
\end{eqnarray}
But $\lambda R^2 = \alpha A^{2/d}$, and it follows from Lemma \ref
{ro1=cs} that
$\forall\alpha>0$,
$\forall A >0$, for $\forall\varepsilon>0$,
%
\begin{eqnarray}
&&\limsup_{T \rightarrow\infty} \frac{1}{T^{1/q}}
\log P\bigl[{\|{Z}\|_{2q,R} \geq\sqrt{2} T^{1/2q} \bigl((m+y+\varepsilon)^{1/2q}
- \varepsilon^{1/2q}\bigr)}\bigr]
\nonumber\\[-8pt]\\[-8pt]
&&\qquad\leq- c(q) \min(1,\alpha A^{2/d})
\bigl((y+\varepsilon)^{1/2q} - \varepsilon^{1/2q} -C(q) A^{1/2q}
\bigr)_+^2.\nonumber
\end{eqnarray}

For the denominator in (\ref{LDUB}), using (\ref{BIZ4.eq}), (\ref{r1=inv_r2})
and part 1 of Lemma \ref{expnablaZ.lem}, we get that
%
\begin{equation}
\liminf_{T \rightarrow\infty} \frac{1}{T^{1/q}}
\log P\bigl[{\|{Z}\|_{2q,R} \geq2 \sqrt{2} T^{1/2q} \varepsilon^{1/2q}}\bigr]
\geq- C(q) \varepsilon^{1/q}.
\end{equation}

We have thus proved that $\forall\alpha>0$,
$\forall A >0$, for $\forall\varepsilon>0$,
%
\begin{eqnarray}
&&\limsup_{T \rightarrow\infty} \frac{1}{T^{1/q}} P[{I_T \geq T y}]
\nonumber\\
&&\qquad\leq\alpha+ C(q) \varepsilon^{1/q} - c(q) \min(1,\alpha A^{2/d})
\\
&&\qquad\quad\hspace*{16.7pt}{}\times
\bigl((y+\varepsilon)^{1/2q} - \varepsilon^{1/2q} -C(q) A^{1/2q} \bigr)_+^2
.\nonumber
\end{eqnarray}

We send $ \varepsilon$ to zero and take $\alpha=A^{-2/d}$, to obtain that
$\forall A >0$,
%
\begin{equation}
\limsup_{T \rightarrow\infty} \frac{1}{T^{1/q}} P[{I_T \geq T y}]
\leq A^{-2/d} - c(q) \bigl(y^{1/2q} -C(q) A^{1/2q} \bigr)_+^2.
\end{equation}
We now choose $A$ such that $C(q) A^{1/2q}= \frac{1}{2} y^{1/2q}$.
$\forall y >0$,
%
\begin{equation}
\limsup_{T \rightarrow\infty} \frac{1}{T^{1/q}} P[{I_T \geq T y}]
\leq- c(q) (y^{1/q} - y^{-2/d} ) \leq-c(q) y^{1/q}
\end{equation}
for $y^{-2/d} \leq y^{1/q} /2$, that is, $y >2$.

\printaddresses

\end{document}